\documentclass[12pt,a4paper,reqno]{amsart}

\usepackage{xcolor}
\usepackage[all,2cell,arrow,color]{xy}
\usepackage[english]{babel}
\usepackage{stmaryrd}
\usepackage{tensor}
\usepackage{graphicx} 
\usepackage{rotating}
\usepackage{pdflscape}
\usepackage[pagebackref,colorlinks,citecolor=blue,linkcolor=blue]{hyperref}

\renewcommand{\wedge}{\circ}

  \theoremstyle{definition}

\makeatletter
\newcommand*{\inlineequation}[2][]{%
  \begingroup
    \refstepcounter{equation}%
    \ifx\\#1\\%
    \else
      \label{#1}%
    \fi
    \relpenalty=10000 %
    \binoppenalty=10000 %
    \ensuremath{%
      #2%
    }%
    ~\@eqnnum
  \endgroup
}

\begin{document}

\title[Hopf polyads, Hopf categories, Hopf group algebras 
as Hopf monads]{Hopf polyads, Hopf categories and Hopf group monoids  
viewed as Hopf monads}

\author{Gabriella B\"ohm} 
\address{Wigner Research Centre for Physics, H-1525 Budapest 114,
P.O.B.\ 49, Hungary}
\email{bohm.gabriella@wigner.mta.hu}

\begin{abstract}
We associate, in a functorial way, a monoidal bicategory $\mathsf{Span}\vert
\mathcal V$ to any monoidal bicategory $\mathcal V$. Two examples of this
construction are of particular interest: {\em Hopf polyads} of \cite{Bru} can
be seen as Hopf monads in $\mathsf{Span}\vert \mathsf{Cat}$ while 
{\em Hopf group monoids} in the spirit of \cite{Z,T} in a braided monoidal 
category $V$, and {\em Hopf categories} of \cite{BCV} over $V$ both 
turn out to be Hopf monads in $\mathsf{Span}\vert V$. Hopf group 
monoids and Hopf categories are Hopf monads on a distinguished type 
of monoidales fitting the framework of \cite{BoLa}. These examples 
are related by a monoidal pseudofunctor $V\to \mathsf{Cat}$.  
\end{abstract}
  
\maketitle


\section{Introduction} \label{sec:intro}

A {\em Hopf monad} \cite{CLS} in a monoidal bicategory is an opmonoidal monad
on a monoidale (also called a pseudo monoid) such that certain fusion 2-cells
are invertible (cf. Section \ref{sec:OpMonoidalMonad}). In the
monoidal 2-category $\mathsf{Cat}$ of categories, functors and natural
transformations, the Hopf monads of \cite{BLV} on monoidal categories are
re-obtained. Opmonoidal monads (in any bicategory) have the
characteristic feature that their Eilenberg-Moore object --- provided that it
exists --- is a monoidale too such that the forgetful morphism is a strict
morphism of monoidales. If the base monoidale is also closed, then the Hopf
property is equivalent to the lifting of the closed structure to the
Eilenberg-Moore object, see \cite{CLS}.   

A monoidale is said to be a {\em map} monoidale if its multiplication and 
unit 1-cells possess right adjoints. We say that it is an {\em opmap} 
monoidale if it is a map monoidale in the vertically opposite 
bicategory (that is, in the original bicategory the multiplication 
and the unit are right adjoints themselves). Thus passing to the 
vertically opposite bicategory, opmonoidal monads on opmap monoidales 
can be seen as monoidal comonads on map monoidales, the central 
objects of the study in \cite{BoLa}. 

An (op)map monoidale is said to be {\em naturally Frobenius}
\cite{Lop-Fr,Lop-FrStreetWood} if two canonical 2-cells (explicitly recalled
in \cite[Paragraph 2.4]{BoLa}), relating the multiplication and its adjoint,
are invertible. The endohom category of a naturally Frobenius (op)map
monoidale in any  monoidal bicategory admits a duoidal structure \cite{St}
(what was called a 2-monoidal structure in \cite{AgMa}). The Hopf monads on a 
naturally Frobenius opmap monoidale can be regarded as Hopf monoids in 
this duoidal endohom category. In this setting, many equivalent 
characterizations --- including the existence of an antipode --- of 
Hopf monads were obtained in \cite{BoLa}. 

Hopf monads in monoidal bicategories unify various structures like groupoids,
Hopf algebras, weak Hopf algebras \cite{BNSz}, Hopf algebroids 
\cite{Sch}, Hopf monads of \cite{BV} and --- more generally ---
of \cite{BLV}. Some of these, namely groupoids, Hopf algebras, weak Hopf
algebras \cite{BNSz}, Hopf algebroids over commutative algebras as in
\cite{Rav} and the Hopf monads of \cite{BV} live on naturally Frobenius opmap 
monoidales, see \cite{BoLa}.

The aim of this note is to show that some structures that recently appeared in
the literature fit this framework as well: we show that {\em Hopf 
group monoids} (thus in particular {\em Hopf group algebras} in 
\cite{T,Z,CL}), {\em Hopf categories} in \cite{BCV} and {\em 
Hopf polyads} in \cite{Bru} can be seen as Hopf monads in suitable 
monoidal bicategories. Hopf group monoids and Hopf categories are 
even Hopf monads on naturally Frobenius opmap monoidales; explaining 
e.g. the existence and the properties of their antipodes. 

Note that all of Hopf polyads, Hopf group monoids, and Hopf categories 
can be seen as lax functors from a suitable category (provided by 
an arbitrary category, a group, and an indiscrete category, respectively) 
to a monoidal bicategory $\mathcal V$ (equal to $\mathsf{Cat}$ and 
a braided monoidal category regarded as a monoidal bicategory with a single
object, respectively); so they are objects of a bicategory of lax functors, 
lax natural transformations and modifications. However, this 
bicategory does not admit a suitable monoidal structure allowing for 
a study of Hopf monads. 

So in order to achieve our goal, we embed it into a larger bicategory 
$\mathsf{Span}\vert \mathcal V$. The bicategory $\mathsf{Span}\vert
\mathcal V$ is constructed for any bicategory $\mathcal V$. 
Whenever $\mathcal V$ is a monoidal bicategory, also $\mathsf{Span}\vert
\mathcal V$ is proven to be so. This correspondence is functorial in the sense
that any lax functor (respectively, monoidal lax functor) $F: \mathcal V \to
\mathcal W$ induces a lax functor (respectively, monoidal lax functor) 
$\mathsf{Span}\vert F: \mathsf{Span}\vert \mathcal V \to \mathsf{Span}\vert
\mathcal W$. This construction is applied to two examples:
\begin{itemize}
\item[{---}]
A monad in $\mathsf{Span}\vert \mathsf{Cat}$ is precisely a {\em polyad} of 
\cite{Bru}.
Furthermore, any set of monoidal categories can be regarded as a monoidale in
$\mathsf{Span}\vert \mathsf{Cat}$.  
The opmonoidal structures of a monad on such a monoidale correspond
bijectively to opmonoidal structures of the polyad in the sense of \cite{Bru}. 
Finally, such an opmonoidal monad is a Hopf monad if and only if the 
corresponding opmonoidal polyad is a Hopf polyad (in the sense of 
\cite{Bru}) over a groupoid. 
\item[{---}]
Any braided monoidal category $V$ can be regarded as a monoidal bicategory
with a single object. Hence there is an associated monoidal bicategory
$\mathsf{Span}\vert V$ in which any  object carries the structure of a
naturally Frobenius opmap monoidale. 

On the one hand, we identify categories enriched in $V$ with certain 
monads; categories enriched in the category of comonoids in $V$ with 
certain opmonoidal monads; and Hopf categories over $V$ with certain 
Hopf monads on these naturally Frobenius opmap monoidales in 
$\mathsf{Span}\vert V$.    

On the other hand, we also identify monoids in $V$ graded by ordinary 
monoids with monads; semi-Hopf group monoids in $V$ with opmonoidal 
monads; and Hopf group monoids in $V$ with Hopf monads on a 
trivial naturally Frobenius opmap monoidale in $\mathsf{Span}\vert V$.    
\end{itemize}
The above examples are related by a monoidal pseudofunctor $V\to 
\mathsf{Cat}$. It induces a monoidal pseudofunctor 
$\mathsf{Span}\vert V \to \mathsf{Span}\vert \mathsf{Cat}$ which 
takes both Hopf group monoids and Hopf categories to Hopf polyads.  

\subsection*{Acknowledgement} I would like to thank Joost Vercruysse for
highly inspiring discussions on the topic and the organizers of the
conference ``StefFest'' in May 2016 in Turin, where this exchange of ideas
began.  
Special thanks go to the anonymous referee for a careful reading of the
manuscript and for helpfully constructive comments. 
Financial support by the Hungarian Scientific Research Fund OTKA (grant
K108384) is gratefully acknowledged.  

\section{The general construction}

Throughout this section $\mathcal V$ will denote a bicategory \cite[Vol. 1
Section 7.7]{Borceaux} whose horizontal composition will be denoted by $\circ$
and whose vertical composition will be denoted by $\ast$. Although the
horizontal composition is required to be neither strictly associative nor
strictly unital, we will omit explicitly denoting the associativity and
unitality iso 2-cells. 

\subsection{Hopf monads in monoidal bicategories}
\label{sec:OpMonoidalMonad}
We briefly recall some definitions for later reference. For more details we 
refer e.g. to \cite{CLS}.

A {\em monad} on a {\em category} $A$ consists of an endofunctor $f:A\to A$ 
together with natural transformations $\mu$ (the multiplication) from 
the two-fold iterate $f\circ f$ to $f$ and $\eta$ (the unit) from the identity 
functor $1$ to $f$. They are subject to the associativity and unitality
axioms.

From the 2-category $\mathsf{Cat}$ of categories, functors and 
natural transformations, this notion can be generalized to {\em any 
bicategory}, see \cite{Str:FTMI}. Then a monad consists of a 1-cell 
$f:A\to A$ and 2-cells $\mu:f\circ f \to f$ and $\eta:1\to f$ such 
that $\mu$ is associative with unit $\eta$.

For {\em monoidal categories} $A$ and $A'$ (with respective monoidal 
products $\otimes$ and $\otimes'$; monoidal units $K$ and $K'$),
we can ask about the relation of a functor $f:A \to A'$ and the 
monoidal structures; there are some dual possibilities of their
compatibility. An {\em opmonoidal} (by some authors called {\em 
comonoidal}) structure on $f$ consists of natural transformations 
$f_2:f(-\otimes  -) \to f(-) \otimes' f(-)$ and $f_0:f(K)\to K'$ 
which satisfy the evident coassociativity and counitality conditions (see
these conditions spelled out explicitly in a more general case below). 
The functor $f$ is said to be {\em strict monoidal} if $f_2$ and $f_0$ are
identity morphisms.

A natural transformation between opmonoidal functors $f$ 
and $f'$ is said to be {\em opmonoidal} if compatible with the 
opmonoidal structures of $f$ and $f'$ (for the precise form of this 
compatibility see the more general case below).

It is straightforward to see that monoidal categories, opmonoidal 
functors and opmonoidal natural transformations constitute a 
2-category $\mathsf{OpMon}$. The monads therein are termed {\em opmonoidal
monads}. Recall from \cite{Moerdijk} and \cite{McCr} that for any monoidal
category $A$ and any monad $f$ on the category $A$, there is a bijective
correspondence between   
\begin{itemize}
\item[{---}]
opmonoidal structures of the functor $f$ making it an opmonoidal monad; 
\item[{---}]
monoidal structures of the category $A^f$ of Eilenberg--Moore $f$-algebras
such that the forgetful functor $A^f \to A$ is strict monoidal 
(that is, the {\em liftings} of the monoidal structure of $A$ to $A^f$).
\end{itemize} 

To any opmonoidal monad $(f,f_2,f_0,\mu,\eta)$ on a monoidal category $A$, one
associates a natural transformation, the so-called {\em fusion morphism},
$$
\xymatrix@C=35pt{
f(f(-) \otimes -) \ar[r]^-{f_2} &
f(f(-)) \otimes f(-) \ar[r]^-{\mu \otimes 1} &
f(-) \otimes f(-).
}
$$
The opmonoidal monad $f$ is said to be a {\em Hopf monad} precisely if the
fusion morphism is invertible, see \cite{BLV}. 
Whenever the monoidal category $A$ is {\em closed}, the invertibility of the
fusion morphism is equivalent to the lifting of the closed structure of $A$ to
the Eilenberg--Moore category $A^f$, see again \cite{BLV}.

The above notions can be generalized from the Cartesian monoidal 
2-category $\mathsf{Cat}$ to any {\em monoidal bicategory} $\mathcal 
V$ (with monoidal product $\otimes$ and monoidal unit $K$). Then monoidal
category is generalized to what is known as {\em monoidale} (alternatively
called {\em pseudo monoid}). Such a gadget consists of an object $A$ of
$\mathcal V$ together with 1-cells $m$ from the monoidal square $A\otimes A$
to $A$ and $u$ from the monoidal unit $K$ to $A$; as well as invertible 
2-cells $m\circ (m\otimes 1) \to m\circ (1\otimes m)$, $m\circ 
(u\otimes 1) \to 1$ and $m\circ (1\otimes u) \to 1$ which satisfy 
McLane's coherence axioms.

For monoidales $A$ and $A'$, an {\em opmonoidal 1-cell} consists of a 
1-cell $f:A\to A'$ together with 2-cells $f_2:f\circ m \to m'\circ 
(f\otimes f)$ and $f_0: f\circ u \to u'$ satisfying the usual 
coassociativity and counitality conditions
$$
\xymatrix@C=40pt@R=20pt{
f\circ m\circ (m\otimes 1) \ar[r]^-{f_2 \circ 1} 
\ar[d]_-{\cong} &
\stackrel{\displaystyle m'\circ (f\otimes f) \circ (m\otimes 1)}
{\cong m'\circ (f\circ m\otimes f)}
\ar[r]^-{1\circ (f_2\otimes 1)} &
\stackrel{\displaystyle m'\circ (m'\otimes 1) \circ (f\otimes f \otimes f) }
{\cong m'\circ (m'\circ (f\otimes f)\otimes f)}
\ar[d]^-\cong \\
f\circ m\circ (1\otimes m) \ar[r]_-{f_2 \circ 1} &
\stackrel{\displaystyle m'\circ (f\otimes f) \circ (1\otimes m)}
{\cong m'\circ (f\otimes f\circ m)}
\ar[r]_-{1\circ (1 \otimes f_2)} &
\stackrel{\displaystyle m'\circ (1\otimes m') \circ (f\otimes f \otimes f) }
{\cong m'\circ (f\otimes m'\circ (f\otimes f))}
\\
f\circ m\circ (u\otimes 1) \ar[r]^-{f_2 \circ 1} 
\ar[d]_-{\cong} &
\stackrel{\displaystyle m'\circ (f\otimes f) \circ (u\otimes 1)}
{\cong m'\circ (f\circ u\otimes f)}
\ar[r]^-{1\circ (f_0\otimes 1)} &
\stackrel{\displaystyle m'\circ (u'\otimes 1) \circ f}
{\cong m'\circ (u'\otimes f)}
\ar[d]^-\cong \\
f \ar@{=}[rr] &&
f \\
f\circ m\circ (1\otimes u) \ar[r]_-{f_2 \circ 1} 
\ar[u]^-{\cong} &
\stackrel{\displaystyle m'\circ (f\otimes f) \circ (1\otimes u)}
{\cong m'\circ (f\otimes f\circ u)}
\ar[r]_-{1\circ (1\otimes f_0)} &
\stackrel{\displaystyle m'\circ (1\otimes u') \circ f}
{\cong m'\circ (f\otimes u').}
\ar[u]_-\cong 
}
$$
A {\em strict monoidal} 1-cell is an opmonoidal 1-cell $f$ with $f_2$ and
$f_0$ the identity 2-cells.

A 2-cell $\varphi:f\to f'$ between opmonoidal 1-cells is {\em opmonoidal}
if the diagrams
$$
\xymatrix{
f\circ m \ar[r]^-{f_2} \ar[d]_-{\varphi \circ 1} &
m'\circ (f\otimes f) \ar[d]^-{1\circ (\varphi \otimes \varphi)} \\
f' \circ m \ar[r]_-{f'_2} &
m' \circ (f' \otimes f')}\qquad
\xymatrix{
f\circ u \ar[r]^-{f_0} \ar[d]_-{\varphi \circ 1} &
u'\ar@{=}[d] \\
f' \circ u \ar[r]_-{f'_0} &
u'}
$$
commute.

Once again, monoidales, opmonoidal 1-cells and opmonoidal 2-cells 
constitute a bicategory $\mathsf{OpMon(\mathcal V)}$; the monads therein are
termed {\em opmonoidal monads}. Assume that in $\mathcal V$ the
Eilenberg-Moore object $A^f$ exists for any monad $f$ on some object $A$.
Then for any monad $f$ on $A$, and for any monoidale with
object part $A$, there is a bijective correspondence between  
\begin{itemize}
\item[{---}] 2-cells $f\circ m \to m\circ (f\otimes f)$ and $f\circ u \to u$
  yielding an opmonoidal monad $f$;  
\item[{---}] 1-cells $A^f\otimes A^f \to A^f$ and $K\to A^f$ yielding a
  monoidale $A^f$ such that the forgetful 1-cell $A^f \to A$ is strict
  monoidal. 
\end{itemize}

The {\em fusion 2-cell} associated to an opmonoidal monad
$(f,f_2,f_0,\mu,\eta)$ takes now the form
$$
\xymatrix@C=35pt{
f\circ m \circ (f\otimes 1)
\ar[r]^-{f_2\circ 1} &
m\circ (f\otimes f) \circ (f\otimes 1) \cong 
m\circ (f\circ f \otimes f)
\ar[r]^-{1\circ (\mu \otimes 1)} &
m\circ (f\otimes f).}
$$
Its invertibility defines $f$ to be a {\em Hopf monad}. As shown in
\cite{CLS}, in the case when the base monoidale is {\em closed}, the
invertibility of the fusion 2-cell is again equivalent to the lifting of the
closed structure to the Eilenberg-Moore object of $f$.
For some equivalent characterizations of Hopf monads (among opmonoidal monads)
in favorable situations, we refer to \cite{BoLa}.

\subsection{The bicategory $\mathsf{Span}\vert \mathcal V$ associated to a
  bicategory $\mathcal V$} 
The {\bf 0-cells} of $\mathsf{Span}\vert \mathcal V$ are pairs 
consisting of a set $X$ and a map $x$ from $X$ to the set $\mathcal 
V^0$ of $0$-cells in $\mathcal V$. 

The {\bf 1-cells} from 
$\xymatrix@C=16pt{
X \ar[r]|(.42){\,x\,} & \mathcal V^0}$ 
to 
$\xymatrix@C=16pt{
Y \ar[r]|(.42){\,y\,} & \mathcal V^0 }$
consist of a span 
$\xymatrix@C=16pt{
Y & A \ar[l]|(.42){\,l\,} \ar[r]|(.42){\,r\,} & X }$ --- inducing a span 
$\xymatrix@C=20pt{
\mathcal V^0 & A \ar[l]|(.42){\,y.l\,} \ar[r]|(.42){\,x.r\,} & \mathcal V^0 }$
--- and a map $a$ from $A$ to the set $\mathcal V^1$ of $1$-cells in 
$\mathcal V$, such that with the source and target maps $s$ and $t$ 
of $\mathcal V$ the following compatibility diagram commutes 
(that is to say, $a$ is a map of spans over the set $\mathcal V^0$).   
\begin{equation}\label{eq:1-cell}
\xymatrix{
Y \ar[d]_-y &
A \ar[l]_-l \ar[r]^-r \ar[d]^-a &
X \ar[d]^-x \\
\mathcal V^0 &
\mathcal V^1 \ar[l]^-t \ar[r]_- s &
\mathcal V^0}
\end{equation}

The {\bf 2-cells} from 
$(\xymatrix@C=12pt{
Y & A \ar[l] \ar[r] & X,}a)$
to 
$(\xymatrix@C=12pt{
Y & A' \ar[l] \ar[r] & X,}a')
$
consist of a map of spans $f:A\to A'$ and a set 
$\varphi=\{\varphi_c: a(c) \Rightarrow a'f(c)\vert c\in A\}$ of 
2-cells in $\mathcal V$. 

If we regard the maps $a$ and $a'$ as functors from the discrete 
categories $A$ and $A'$, respectively, to the vertical category of $\mathcal
V$, then $\varphi$ is a natural transformation from $a$ to the composite 
of the functors $f:A \to A'$ and $a'$. By this motivation we  
use the diagrammatic notation
$$
\xymatrix{
A \ar[rr]^-a \ar[rd]_-f &
\ar@{}[d]|-{\Downarrow\varphi} &
\mathcal V^1. \\
&
A' \ar[ru]_{a'}
}
$$

The {\bf vertical composite} of the 2-cells $(f,\varphi):
(\xymatrix@C=10pt{
Y & A \ar[l]\ar[r] & X,}a) 
\Rightarrow  
(\xymatrix@C=10pt{
Y & A' \ar[l] \ar[r] & X,}$ $a')$
and  
$(f',\varphi'):
(\xymatrix@C=10pt{
Y & A' \ar[l] \ar[r] & X,}a')
\Rightarrow
(\xymatrix@C=10pt{
Y & A'' \ar[l] \ar[r] & X,}a'')$ 
is the pair
$$
\xymatrix{
& A \ar[ld] \ar[rd] \ar[d]^-f\\
Y & A'\ar[l] \ar[d]^-{f'} \ar[r] & X \\
& A'' \ar[lu] \ar[ru]}
\qquad
\xymatrix@C=65pt@R=65pt{
A \ar[r]^-a \ar[d]_f 
\ar@{}[rd]|(.25){\Downarrow \varphi}|(.75){\Downarrow \varphi'}&
\mathcal V^1 \\
A' \ar[ru]^-{a'} \ar[r]_{f'} &
A''. \ar[u]_-{a''}}
$$
In other words, it is the pair $(f'.f, \{\varphi'_{f(c)} \ast 
\varphi_c \vert c \in A\})$. 

The identity 2-cell of 
$(\xymatrix@C=12pt{
Y & A \ar[l] \ar[r] & X,}a)$
consists of the identity map $1:A \to A$ and the set $\{ 1_{a(c)} 
\vert c\in A\}$ of identity 2-cells.

The {\bf horizontal composite} of the 1-cells 
$(\xymatrix@C=15pt{
Y & A \ar[l]|(.4){\, l \,} \ar[r]|(.4){\, r \,} & X,}a)$
and 
$(\xymatrix@C=15pt{
Z & B \ar[l]|(.4){\, l \,} \ar[r]|(.4){\, r \,} & Y,}$ $b)$ 
is the pair consisting of the pullback span 
$$
Z \leftarrow B\wedge A :=\{(d,c)\in B\times A \vert r(d)=l(c)\} \to 
X, \qquad 
l(d) \mapsfrom (d,c) \mapsto r(c) 
$$
and the map 
$$
B\wedge A\to \mathcal V^1,\qquad (d,c) \mapsto b(d)\circ a(c). 
$$
The 1-cells $b(d)$ and $a(c)$ are composable indeed thanks to 
(\ref{eq:1-cell}).

The horizontal composite of 2-cells 
$(f,\varphi): 
(\xymatrix@C=12pt{
Y & A \ar[l] \ar[r] & X,}a)
\Rightarrow 
(\xymatrix@C=12pt{
Y & A' \ar[l] \ar[r] & X,}a')
$ and 
$(g,\gamma): 
(\xymatrix@C=12pt{
Z & B \ar[l] \ar[r] & Y,}b) 
\Rightarrow 
(\xymatrix@C=12pt{
Z & B' \ar[l] \ar[r] & Y,}b')$ consists of the map 
$$
g \wedge f: B \wedge A \to B'\wedge A',\qquad
(d,c) \mapsto (g(d),f(c))
$$
and the following set of 2-cells in $\mathcal V$.
$$
\{\gamma_d \circ \varphi_c: b(d)\circ a(c) \Rightarrow b'g(d) \circ 
a'f(c) \vert (d,c) \in B \wedge A\}
$$
The identity 1-cell of $(X,x)$ consists of the trivial span 
$\xymatrix@C=12pt{
X \ar@{=}[r] &
X \ar@{=}[r] &
X}$ 
and the map $1_{x(-)}:X \to \mathcal V^1$. 
The associativity and unitality natural transformations are pairs of the 
analogous natural transformations in $\mathsf{Span}$ and $\mathcal V$.

Using that both $\mathsf{Span}$ and $\mathcal V$ are bicategories, it 
is straightforward to see that so is $\mathsf{Span}\vert \mathcal V$ 
above.

We are not aware of any construction yielding $\mathsf{Span}\vert 
\mathcal V$ as a comma bicategory. However, regarding it as a 
tricategory (with only identity 3-cells), it embeds into a comma 
tricategory obtained by a lax version of the {\em 3-comma category 
construction} in \cite[Section I.2.7]{Gray}: 
Consider the tricategory $\mathsf{SpanSpan}$ whose {0-cells} are sets
$X,Y,\dots$,
whose hom-bicategory $\mathsf{SpanSpan}(X,Y)$ is the bicategory of spans
in the category $\mathsf{Span}(X,Y)$, and in which the 
{1-composition} is the pullback of spans with the evident coherence 
2-and 3-cells. Regarding $\mathsf{Span}$ as a tricategory with only 
identity 3-cells, and interpreting a map of spans in the first diagram of
$$
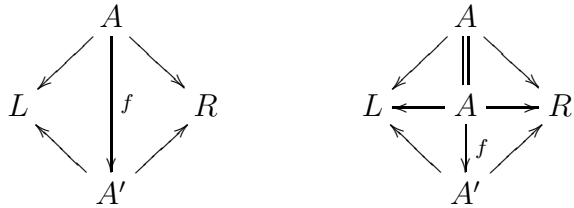
\xymatrix @C=18pt @R=18pt {
&A \ar[ld] \ar[rd] \ar[dd]^-f\\
L & & R \\
& A' \ar[lu] \ar[ru]}
\qquad \qquad \xymatrix @C=18pt @R=18pt{
&A \ar[ld] \ar[rd] \ar@{=}[d]\\
L & A \ar[l] \ar[r] \ar[d]^-f & R \\
& A' \ar[lu] \ar[ru]}
$$
as a span in the second diagram, we obtain a functor of tricategories 
$\mathsf{Span} \to \mathsf{SpanSpan}$. 
On the other hand, any bicategory $\mathcal V$ determines an evident 
(1- and 2-) lax functor of tricategories from the trivial 
tricategory $\mathbf{1}$ (with a single 0-cell and only identity
higher cells) to $\mathsf{SpanSpan}$.   
The comma tricategory arising from the lax functors 
$\xymatrix@C=16pt{
\mathsf{Span} \ar[r] &
\mathsf{SpanSpan} &
\mathbf{1} \ar[l]|(.26){\,{\mathcal{V}} \,}}$
contains $\mathsf{Span}\vert \mathcal V$ as a sub-tricategory. 

Note for later application that a 1-cell $(\xymatrix@C=12pt{Y & A \ar[l]
\ar[r] & X,}a)$ possesses a right adjoint in $\mathsf{Span}\vert \mathcal V$
if and only if $\xymatrix@C=12pt{Y & A \ar[l] \ar[r] & X}$ has a right adjoint
in $\mathsf{Span}$ and for all $c\in A$, $a(c)$ has a right adjoint in $\mathcal
V$. Equivalently, if and only if it is isomorphic to a 1-cell of the form 
$(\xymatrix@C=12pt{Y & X \ar[l] \ar@{=}[r] & X,}h)$ such that for all $p\in
X$, $h(p)$ has a right adjoint in $\mathcal V$. 

\subsection{Monads in $\mathsf{Span}\vert {\mathcal V}$} \label{sec:monad}
Let us fix an arbitrary 0-cell 
$(D^0, 
\xymatrix@C=15pt{D^0 \ar[r]|(.45)f & \mathcal V^0})$ in 
$\mathsf{Span}\vert \mathcal V$ and describe a monad on it. 
The underlying 1-cell consists of a span 
$\xymatrix@C=16pt{
D^0 & D^1 \ar[l]|(.45){\, t \,} \ar[r]|(.45){\, s \,} & D^0}$
and a map $F$ associating a 1-cell $F(h):fs(h) \to ft(h)$ in $\mathcal V$ to
each element $h$ of $D^1$. The multiplication and unit 2-cells consist of 
respective maps of spans
$$
\xymatrix@C=10pt{
& D^1 \wedge D^1 \ar[ld] \ar[rd] \ar[dd]^-{\cdot}\\
D^0 && D^0 \\
& D^1 \ar[lu] \ar[ru]}
\qquad 
\raisebox{-37pt}{$\textrm{and}$}
\qquad
\xymatrix{
& D^0 \ar@{=}[ld] \ar@{=}[rd] \ar[dd]^-{e}\\
D^0 && D^0 \\
& D^1 \ar[lu] \ar[ru]}
$$
and respective sets of 2-cells $\{ \mu_{h,k}: 
F(h)\circ F(k) \to F(h.k) \vert (h,k) \in D^1\wedge D^1\}$ and $\{ \eta_x: 
1_{f(x)} \to F(e_x) \vert x \in D^0\}$ in $\mathcal V$. The associativity and
unitality conditions precisely say that there is a category 
\begin{equation}\label{eq:D-cat}
\xymatrix{
D^0 \ar[r]|-{\,e\,} &
D^1 \ar@<-5pt>[l]_-s \ar@<5pt>[l]^-t &
D^1 \circ D^1 \ar[l]_(.55){\cdot}}
\end{equation}
with object set $D^0$, morphism 
set $D^1$, source and target maps $s$ and $t$, composition~$\cdot$ 
and identity morphisms $\{e_x\vert x\in D^0\}$ and --- regarding this category
as a bicategory with only identity 2-cells --- a lax functor $D\to \mathcal V$
with object map $f$, hom functor $F$ (from the discrete hom category $D^1$),
and comparison natural transformations $\mu$ and $\eta$. Summarizing, for any
bicategory $\mathcal V$, the following notions coincide. 
\begin{itemize}
\item[{---}] A pair consisting of a category $D$ and a lax functor $D \to
  \mathcal V$.
\item[{---}] A monad in $\mathsf{Span}\vert {\mathcal V}$.
\end{itemize}

\subsection{Bicategories of monads in $\mathsf{Span}\vert {\mathcal
    V}$} \label{sec:lax_funtor_bicat}  
Consider a category (\ref{eq:D-cat}) and lax functors $((f,F),\mu,\eta)$ and
$((f',F'),\mu',\eta')$ from $D$ to ${\mathcal V}$. Regard them as monads in
$\mathsf{Span}\vert {\mathcal V}$ as in Section \ref{sec:monad}.

A 1-cell of the form
$(\xymatrix@C=10pt{D^0 & D^0 \ar@{=}[l] \ar@{=}[r] & D^0}, 
  \xymatrix@C=18pt{D^0 \ar[r]|(.5){\,h\,} & \mathcal V^1})$
from 
$\xymatrix@C=18pt{
D^0 \ar[r]|(.5){\,f\,} & \mathcal V^0}$ 
to 
$\xymatrix@C=18pt{
D^0 \ar[r]|(.45){\,f'\,} & \mathcal V^0}$ and
the 2-cell 
$(D^0\circ D^1 \cong D^1\cong D^1\circ D^0,
\xymatrix@C=18pt{ht(-) \circ F(-) \ar[r]|(.47){\,\varphi\,} & 
F'(-) \circ hs(-)})$
constitute a monad morphism (in the sense of \cite{Str:FTMI}) in
$\mathsf{Span}\vert {\mathcal V}$ if and only if $(h,\varphi)$ is a lax
natural transformation. 

A 2-cell of the form
$(D^0=D^0,
\xymatrix@C=18pt{h\ar[r]|(.45){\,\gamma\,} & h'})$ is a monad transformation
(in the sense of \cite{Str:FTMI}) in $\mathsf{Span}\vert {\mathcal V}$ if and
only if $\gamma$ is a modification $(h,\varphi) \to (h',\varphi')$.

These observations amount to the isomorphism of the following bicategories,
for any category $D$ and any bicategory $\mathcal V$.
\begin{itemize}
\item[{---}] The bicategory $[D,\mathcal V]$ of lax functors $D \to \mathcal
  V$, lax natural transformations and modifications.
\item[{---}] The following locally full sub-bicategory
in the bicategory of monads in $\mathsf{Span}\vert \mathcal V$. 
The 0-cells are those monads which live on 0-cells $D^0 \to \mathcal V^0$ 
(for the given object set $D^0$ of $D$), whose 1-cells are of the form 
$(\xymatrix@C=14pt{
\! D^0 & D^1 \ar[l]|(.45){\, t \,} \ar[r]|(.45){\, s \,} & D^0,}$ $F)$ 
(in terms of the given data $s,t$), and whose multiplication and unit 
2-cells have the respective forms $(\cdot,\mu)$ and $(e,\eta)$ (with 
the given maps $\cdot$ and $e$). The 1-cells are those monad morphisms
$((H,h),(f,\varphi))$ whose underlying span $H$ is the trivial  
span $D^0=D^0=D^0$ and whose map $f$ is the canonical isomorphism $D^0\circ
D^1 \cong D^1 \cong D^1 \circ D^0$. The 2-cells are all possible monad
transformations $(g,\gamma)$ ($g$ in them is necessarily the identity map
$D^0\to D^0$).  
\end{itemize}

\subsection{The monoidal bicategory $\mathsf{Span}\vert \mathcal V$ for a
monoidal  bicategory $\mathcal V$} 
In this section $\mathcal V$ is taken to be a monoidal bicategory --- that is,
a single object tricategory \cite{GPS} --- with monoidal operation 
$\otimes$ and monoidal unit $K$. Then we can equip 
$\mathsf{Span}\vert \mathcal V$ with a monoidal structure as follows. 

The monoidal product of 0-cells 
$\xymatrix@C=16pt{
X \ar[r]|(.42){\,x\,} & \mathcal V^0}$ 
and
$\xymatrix@C=16pt{
Z \ar[r]|(.42){\,z\,} & \mathcal V^0 }$
consists of the Cartesian product set $X\times Z$ and the map 
$$
X\times Z \to \mathcal V^0, \qquad (k,l) \mapsto x(k) \otimes z(l).
$$
The monoidal product lax functor on the local hom categories takes a pair of 
2-cells  
$(f,\varphi):
(\xymatrix@C=8pt{
Y & A \ar[l] \ar[r] & X,}a)
\Rightarrow
(\xymatrix@C=8pt{
Y & A' \ar[l] \ar[r] & X,}a')$ 
and 
$(g,\gamma): 
(\xymatrix@C=8pt{
W & B \ar[l] \ar[r] & Z,}b) 
\Rightarrow
(\xymatrix@C=8pt{
W & B' \ar[l] \ar[r] & Z,}b')$
to the 2-cell consisting of the Cartesian product map $f \times g:A \times B 
\to A'\times B'$ between the Cartesian product spans and the following 
set of 2-cells.  
$$
\{ \varphi_c \otimes \gamma_d:a(c) \otimes b(d) \to a'f(c) \otimes 
b'g(d)\vert (c,d)\in A \times B\}
$$
The natural transformations establishing the compatibility of these functors 
$\otimes$ with the identity 1-cells and the horizontal composition are
inherited from $\mathcal V$. 
The monoidal unit is the singleton set with the map $K$.
The lax natural transformations measuring the non-associativity and 
non-unitality of $\otimes$, as well as their invertible coherence modifications 
are induced by those in $\mathcal V$.

It requires some patience to check that this is a monoidal bicategory 
indeed. No conceptual difficulties arise, however, one has to use 
repeatedly that $\mathsf{Span}$ is a monoidal bicategory via the 
Cartesian product of sets together with the assumed monoidal 
bicategory structure of $\mathcal V$.

Note that the sub-bicategory of $\mathsf{Span}\vert \mathcal V$ 
occurring in Section \ref{sec:lax_funtor_bicat} is not a monoidal 
sub-bicategory. Hence it is not suitable for our study of Hopf 
monads. 

\subsection{The bicategory $\mathsf{Span}\vert \mathsf{OpMon}(\mathcal V)$ for
a monoidal bicategory $\mathcal V$} \label{sec:Span|OpMon}
\label{sec:induced_monoidale} 
The 2-full (i.e. both horizontally and vertically full) sub-bicategory in
$\mathsf{OpMon}(\mathsf{Span})$ whose objects are the opmap monoidales, is in
fact isomorphic to $\mathsf{Span}$ via the forgetful functor. 

Consider next a monoidale in $\mathsf{Span}\vert \mathcal V$ whose
multiplication and unit 1-cells have underlying spans which possess left
adjoints in Span. It consists of a 0-cell 
$\xymatrix@C=18pt{X \ar[r]|(.45){\, C \,}& \mathcal V^0}$ 
together with multiplication and unit 1-cells which must be of the form 
\begin{equation}\label{eq:opmap-monoidale}
(\xymatrix@C=12pt{
X & X \ar@{=}[l] \ar[r]^-{\Delta} & X\times X,} 
\xymatrix@C=12pt{X\ar[r]^-m & \mathcal V^1})
\quad 
\textrm{and}
\quad
(\xymatrix@C=12pt{
X & X \ar@{=}[l] \ar[r]^-{!} & 1,} 
\xymatrix@C=12pt{X \ar[r]^-u & \mathcal V^1})
\end{equation}
--- 
where $\Delta$ is the diagonal map $p\mapsto (p,p)$ and $!$ denotes 
the unique map to the singleton set $1$ --- and associativity and 
unit 2-cells provided by the identity map of $X$, and maps sending 
$p\in X$ to 2-cells in $\mathcal V$, $\alpha_p:m_p\circ(m_p \otimes 
1_{C_p}) \to m_p\circ (1_{C_p}\otimes m_p)$, $\lambda_p:m_p\circ (u_p 
\otimes 1_{C_p}) \to 1_{C_p}$ and $\varrho_p:m_p\circ (1_{C_p} 
\otimes u_p) \to 1_{C_p}$, respectively. The axioms for these data to 
constitute a monoidale in $\mathsf{Span}\vert \mathcal V$ say 
precisely that $(C_p,m_p,u_p,\alpha_p,\lambda_p,\varrho_p)$ is a 
monoidale in $\mathcal V$ for all $p\in X$. 

If each member $(C_p,m_p,u_p,\alpha_p,\lambda_p,\varrho_p)$ in a 
monoidale as in (\ref{eq:opmap-monoidale}) is a naturally Frobenius 
opmap monoidale in $\mathcal V$, then so is the induced monoidale in 
$\mathsf{Span}\vert \mathcal V$. The left adjoints of its 
multiplication and unit are given in terms of the left adjoints 
$(m_p)_\ast\dashv m_p$ and $(u_p)_\ast \dashv u_p$ as 
$$ 
(\xymatrix@C=12pt{
X\times X & X \ar[l]_-\Delta \ar@{=}[r]  & X},
X \ni p \mapsto (m_p)_\ast)
\quad
\textrm{and}
\quad
(\xymatrix@C=12pt{
1 & X \ar[l]_-{!} \ar@{=}[r]  & X},
X \ni p \mapsto (u_p)_\ast).
$$
The 2-cells of \cite[Paragraph 2.4]{BoLa} are invertible for the 
induced opmap monoidale since they are so for each member 
$(C_p,m_p,u_p,\alpha_p,\lambda_p,\varrho_p)$. 

In a symmetric manner, a set $\{(C_p,d_p,e_p, \alpha_p,\lambda_p, 
\varrho_p) \vert p\in X\}$ of comonoidales in $\mathcal V$ induces a 
comonoidale in $\mathsf{Span}\vert \mathcal V$. The underlying 0-cell 
is $(X,X\ni p \mapsto C_p)$; the comultiplication and counit 1-cells 
are  
$$ 
(\xymatrix@C=12pt{
X\times X & X \ar[l]_-\Delta \ar@{=}[r]  & X},
X \ni p \mapsto d_p)
\quad
\textrm{and}
\quad
(\xymatrix@C=12pt{
1 & X \ar[l]_-{!} \ar@{=}[r]  & X},
X \ni p \mapsto e_p),
$$
respectively; while the coassociativity and the counit isomorphisms are 
given by the sets $\{ \alpha_p \ \vert \ p\in X\}$, $\{ \lambda_p \ 
\vert \ p\in X\}$ and $\{ \varrho_p \ \vert \ p\in X\}$ of the 
analogous 2-cells for $C_p$.

An opmonoidal 1-cell between monoidales $(X,C)$ and $(Y,H)$ of the form in
(\ref{eq:opmap-monoidale}) consists of a span  
$\xymatrix@C=16pt{
Y & A \ar[l]|(.45){\, l \,} \ar[r]|(.45){\, r \,} & X}$
and a map $a$ sending each element $h$ of $A$ to a 1-cell 
$a(h):C_{r(h)} \to H_{l(h)}$ in $\mathcal V$; together with an opmonoidal
structure which consists of opmonoidal structures on each 1-cell $a(h)$ for
$h\in A$. A 2-cell $(f,\varphi)$ between opmonoidal 1-cells 
$(\xymatrix@C=12pt{
Y & A \ar[l] \ar[r] & X,}a)$
and 
$(\xymatrix@C=12pt{
Y & A' \ar[l] \ar[r] & X,}$ $a')$
as above is opmonoidal precisely if each component $\varphi_h:a(h) 
\to a'f(h)$ is opmonoidal, for $h\in A$. 

Putting in other words, from the considerations of the previous paragraph
isomorphism of the following bicategories follows.
\begin{itemize}
\item[{---}] $\mathsf{Span}\vert \mathsf{OpMon}(\mathcal V)$.
\item[{---}] The 2-full sub-bicategory of $\mathsf{OpMon}(\mathsf{Span}\vert
  \mathcal V)$ whose objects are of the kind in (\ref{eq:opmap-monoidale}).
\end{itemize}

\subsection{Bicategories of monads in $\mathsf{Span}\vert 
\mathsf{OpMon}(\mathcal V)$} \label{sec:Mnd(Span|OpMon)}

Combining the isomorphisms of Section \ref{sec:lax_funtor_bicat} and Section
\ref{sec:Span|OpMon}, we obtain isomorphism of the following
bicategories, for any category $D$ and any monoidal bicategory $\mathcal V$.
\begin{itemize}
\item[{---}] $[D,\mathsf{OpMon}(\mathcal V)]$.
\item[{---}] 
The following locally full sub-bicategory in the bicategory of monads (in the
sense of \cite{Str:FTMI}) in $\mathsf{Span}\vert \mathsf{OpMon}(\mathcal V)$. 
The 0-cells are those monads which live on 0-cells $D^0 \to 
\mathsf{OpMon}(\mathcal V)^0$ 
(for the given object set $D^0$ of $D$), whose 1-cells are of the form 
$(\xymatrix@C=16pt{
D^0 & D^1 \ar[l]|(.45){\, t \,} \ar[r]|(.45){\, s \,} & D^0}, 
\mathbf{d}:D^1 \to \mathsf{OpMon}(\mathcal V)^1)$ 
(in terms of the given data $s,t$), and whose multiplication and unit 
2-cells have the respective forms $(\cdot,\mu)$ and $(e,\eta)$ (with 
the given maps $\cdot$ and $e$).  The 1-cells are those monad
morphisms $((H,\mathbf{h}),(f,\varphi))$ whose underlying span $H$ is the
trivial span $D^0=D^0=D^0$ (hence $\mathbf{h}$ is a map $D^0\to
\mathsf{OpMon}(\mathcal V)^1$) and whose map $f$ is the canonical isomorphism
$D^0\circ D^1 \cong D^1 \cong D^1 \circ D^0$. The 2-cells
are all possible monad transformations $(g,\gamma)$ ($g$ in them is
necessarily the identity map $D^0\to D^0$).  
\item[{---}] The following locally full sub-bicategory in the bicategory of
  monads (in the sense of \cite{Str:FTMI}) in
  $\mathsf{OpMon}(\mathsf{Span}\vert \mathcal V)$.  
The 0-cells are those monads which live on monoidales 
with object part $D^0 \to \mathcal V^0$ 
(for the given object set $D^0$ of $D$) and with multiplication and unit of
the form in (\ref{eq:opmap-monoidale}), whose 1-cells are of the
form  
$(\xymatrix@C=16pt{ 
D^0 & D^1 \ar[l]|(.45){\, t \,} \ar[r]|(.45){\, s \,} & D^0}, d:D^1 \to
\mathcal V^1)$ 
(in terms of the given data $s,t$), and whose multiplication and unit  
2-cells have the respective forms $(\cdot,\mu)$ and $(e,\eta)$ (with 
the given maps $\cdot$ and $e$). (There are no restrictions on the opmonoidal
structure of the 1-cell $(\xymatrix@C=16pt{ 
D^0 & D^1 \ar[l]|(.45){\, t \,} \ar[r]|(.45){\, s \,} & D^0}, d)$
in $\mathsf{Span}\vert \mathcal V$.) The
1-cells are those monad morphisms $((H,h),(f,\varphi))$ whose underlying span
$H$ is the trivial span $D^0=D^0=D^0$ and whose map $f$ is the canonical
isomorphism $D^0\circ D^1 \cong D^1 \cong D^1 \circ D^0$. (There are no
restrictions on the opmonoidal structure of the 1-cell $(D^0=D^0=D^0,h:D^0\to
\mathcal V^1)$ in $\mathsf{Span}\vert \mathcal V$.) The 2-cells are all
possible monad transformations $(g,\gamma)$ ($g$ in them is necessarily the
identity map $D^0\to D^0$).
\end{itemize}

\subsection{Functoriality}
Any lax functor $F:\mathcal V \to \mathcal W$ between arbitrary 
bicategories $\mathcal V$ and $\mathcal W$ induces a lax functor 
$\mathsf{Span} \vert F: \mathsf{Span}\vert \mathcal V \to 
\mathsf{Span} \vert \mathcal W$ as follows. It sends a 0-cell 
$\xymatrix@C=16pt{
X \ar[r]|(.42){\,x\,} & \mathcal V^0}$ 
to the 0-cell
$$
\xymatrix{
X \ar[r]^-x & 
\mathcal V^0 \ar[r]^-{F^0} & 
\mathcal W^0,}
$$ 
and it sends a 2-cell 
$(f,\varphi):
(\xymatrix@C=12pt{
Y & A \ar[l] \ar[r] & X,}a)
\Rightarrow
(\xymatrix@C=12pt{
Y & A' \ar[l] \ar[r] & X,}a')$ 
to 
$$
(f,\{F(\varphi_c) \vert c\in A\}):
(\xymatrix@C=12pt{
Y & A \ar[l] \ar[r] & X,}F(a(-)))
\Rightarrow
(\xymatrix@C=12pt{
Y & A' \ar[l] \ar[r] & X,}F(a'(-))).
$$
The natural transformations establishing its compatibility with the horizontal 
composition and the identity 1-cells come from those for $F$. Hence if $F$ is
a pseudofunctor then so is $\mathsf{Span} \vert F$. 

If $\mathcal V$ and $\mathcal W$ are monoidal bicategories and $F$ is 
a monoidal lax functor (cf. \cite[Definition 3.1]{GPS}) then a monoidal 
structure is induced on $\mathsf{Span} \vert F$ in a natural way.
All the needed axioms hold for $\mathsf{Span}\vert F$ 
thanks to the fact that they hold for $F$. 

Since any lax functor preserves monads, so does $\mathsf{Span}\vert 
F$ for any lax functor $F$. 
Since any monoidal lax functor preserves monoidales, so does 
$\mathsf{Span}\vert F$ for any monoidal lax functor $F$.
Any monoidal lax functor whose unit- and product-compatibilities are 
pseudonatural transformations preserves opmonoidal 1- and 2-cells. 
Hence so does $\mathsf{Span}\vert F$ whenever the unit- and 
product-compatibilities of $F$ are invertible.

\subsection{Convolution monoidal hom categories and their opmonoidal 
monads} \label{sec:conv_mod}

If $M$ is a monoidale and $C$ is a comonoidale in any monoidal 
bicategory $\mathcal M$ then the hom category $\mathcal M(C,M)$ 
admits a monoidal structure of the convolution type: the monoidal 
product of 2-cells $\gamma: b\Rightarrow b'$ and $\varphi: a
\Rightarrow a'$ between 1-cells $C \to M$ is obtained taking the 
horizontal composite of the comultiplication of the 
comonoidale $C$ (which is a 1-cell from $C$ to $C\otimes C$) with the 
monoidal product of $\gamma$ and $\varphi$ in $\mathcal V$ (which 
goes from $C\otimes C$ to $M\otimes M$) and with the multiplication 
of the monoidale $M$ (which is a 1-cell from $M\otimes M$ to $M$). 
The monoidal unit is the horizontal composite of the counit $C\to K$ 
with the unit $K \to M$. 

Via horizontal composition any monad $a: M\to M$ in any bicategory $\mathcal
M$ induces a monad $\mathcal M(C,a)$ in $\mathsf{Cat}$ on the hom category
$\mathcal M(C,M)$, for any 0-cell $C$ of $\mathcal M$. If $C$ is a
comonoidale, $M$ is a monoidale, and $a$ is an opmonoidal monad in $\mathcal
M$, then $\mathcal M(C,a)$ is canonically an opmonoidal monad in 
$\mathsf{Cat}$ on the above convolution-monoidal category $\mathcal 
M(C,M)$. Moreover, if $a$ is a left or right Hopf monad in $\mathcal 
M$ in the sense of \cite{CLS}, then $\mathcal M(C,a)$ is a left or 
right Hopf monad in $\mathsf{Cat}$ in the sense of \cite{BLV}.  

These considerations apply, in particular, to an induced monoidale  
$(Y,M):=\{M_p\vert p\in Y\}$ and an induced comonoidale $(X,C): 
=\{C_q\vert q\in X\}$ in $\mathsf{Span}\vert \mathcal V$ (cf. Section
\ref{sec:induced_monoidale}) for any monoidal bicategory $\mathcal V$. In the
category $\mathsf{Span}\vert \mathcal V((X,C),(Y,M))$ the monoidal product any
two morphisms --- that is, of 2-cells $(g,\gamma): (B,b)\Rightarrow (B',b')$
and $(f,\varphi): (A,a) \Rightarrow (A',a')$ between 1-cells $(X,C) \to (Y,M)$
--- is the morphism consisting of the map of spans 
\begin{equation} \label{eq:prod_span}
\xymatrix{
& 
B\bullet A\colon = 
\{ (c,h) \in B\times A \vert 
l(c)=l(h)\ \textrm{and}\ r(c)=r(h) \} 
\ar[ld]_-{(c,h)\mapsto l(h)} \ar[rd]^-{(c,h)\mapsto r(h)} 
\ar[dd]^-{(c,h)\mapsto (g(c),f(h))}  \\ 
Y && X \\
& B' \bullet A' \ar[lu] \ar[ru]}
\end{equation}
and the set 
$$
\{1 \circ(\gamma_c \otimes \varphi_h)\circ 1 :
m_{l(c)}\circ(b(c)\otimes a(h))\circ d_{r(h)} \to
m_{l(c)}\circ(b'g(c) \otimes  a'f(h)) \circ d_{r(h)} \}
$$ 
of 2-cells in $\mathcal V$ labelled by the elements $(c,h)\in B 
\bullet A$. The monoidal unit $J$ consists of the complete span
$\xymatrix@C=12pt{ 
Y & Y \times X \ar[l] \ar[r] & X}$ 
(whose maps are the first and the second projection, respectively), 
and the map sending $(i,j)\in Y\times X$ to the 1-cell $u_i\circ 
e_j:C_j \to M_i$ in $\mathcal V$.

Now if $(A,a)$ is an opmonoidal monad on $(Y,M)$, then $\mathsf{Span}\vert 
\mathcal V((X,C),(A,a))$ is an opmonoidal monad in $\mathsf{Cat}$ on the 
above monoidal category $\mathsf{Span}\vert \mathcal V((X,C),(Y,M))$; which 
belongs to the realm of the theory of opmonoidal monads in 
\cite{BLV}. 

\section{Hopf polyads as Hopf monads}

In this section we apply the general construction of the previous section to
the 2-category $\mathsf{Cat}$ of categories, functors and natural 
transformations; with the monoidal structure provided by the Cartesian 
product.

\subsection{Monads in $\mathsf{Span}\vert \mathsf{Cat}$ versus polyads}
\label{sec:Span|Cat-monad}
From Section \ref{sec:monad} we conclude on the coincidence of the
following notions.  
\begin{itemize} 
\item[{---}] A {\em polyad} in \cite{Bru}; that is, a pair consisting of a
  category and -- regarding this category as a bicategory with only identity
  2-cells -- a lax functor from it to $\mathsf{Cat}$ (see \cite[Remark
  2.1]{Bru}). 
  \item[{---}] A monad in $\mathsf{Span}\vert \mathsf{Cat}$.
\end{itemize}

By the application of Section \ref{sec:lax_funtor_bicat},
the following bicategories are isomorphic, for any given category
(\ref{eq:D-cat}).  
\begin{itemize}
\item[{---}] The bicategory of polyads over the category
  (\ref{eq:D-cat}) in \cite[Section 3]{Bru}. That is, the bicategory of lax
  functors from (\ref{eq:D-cat}) to $\mathsf{Cat}$, lax natural
  transformations and modifications.
\item[{---}]
The following locally full sub-bicategory in the bicategory of monads (in the
sense of \cite{Str:FTMI}) in $\mathsf{Span}\vert \mathsf{Cat}$. 
The 0-cells are those monads which live on 0-cells $D^0 \to \mathsf{Cat}^0$ 
(for the given object set $D^0$), whose 1-cells are of the form 
$(\xymatrix@C=16pt{
D^0 & D^1 \ar[l]|(.45){\, t \,} \ar[r]|(.45){\, s \,} & D^0}, 
d:D^1\to \mathsf{Cat}^1)$ 
(in terms of the given data $s,t$), and whose multiplication and unit 
2-cells have the respective forms $(\cdot,\mu)$ and $(e,\eta)$ (with 
the given maps $\cdot$ and $e$). The 1-cells are those monad 
morphisms $((H,h),(f,\varphi))$ whose underlying span $H$ is the trivial 
span $D^0=D^0=D^0$ and whose map $f$ is the canonical isomorphism $D^0 
\circ D^1 \cong D^1 \cong  D^1 \circ D^0$. The 2-cells are all possible monad 
transformations $(g,\gamma)$ ($g$ in them is necessarily the identity 
map $D^0\to D^0$).  
\end{itemize}

\subsection{The induced monad in $\mathsf{Cat}$} 
\label{sec:Cat_induced_monad}
Since a polyad is eventually a monad $(D^1,d)$ in $\mathsf{Span}\vert 
\mathsf{Cat}$ on some 0-cell $(D^0,C)$, it induces a monad 
$\mathsf{Span}\vert \mathsf{Cat}((Y,H),(D^1,d))$ in $\mathsf{Cat}$ on 
the category $\mathsf{Span}\vert \mathsf{Cat}((Y,H),(D^0,C))$ for any 
0-cell $(Y,H)$ of $\mathsf{Span}\vert \mathsf{Cat}$, see Section 
\ref{sec:conv_mod}. An object of the Eilenberg-Moore category of this induced 
monad is a pair consisting of a 1-cell $(Q,q):(Y,H) \to (D^0,C)$, and a 
2-cell $(r,\varrho):(D^1,d)\circ (Q,q)\Rightarrow (Q,q)$ in $\mathsf{Span} 
\vert \mathsf{Cat}$ which satisfy the associativity and unitality 
conditions. The morphisms are 2-cells $(Q,q)\Rightarrow (Q',q')$ in 
$\mathsf{Span} \vert \mathsf{Cat}$ which are compatible with the 
actions $(r,\varrho)$ and $(r',\varrho')$.

Let us consider the particular case when the above $Y$ is the 
singleton set $1$ and $H$ takes its single element to the terminal 
category $\mathbf{1}$; and the corresponding Eilenberg-Moore category 
of the monad $\mathsf{Span}\vert \mathsf{Cat}((1,\mathbf{1}),(D^1,d))$.
For any monad $(D^1,d)$ on any 0-cell $(D^0,C)$ in $\mathsf{Span}\vert  
\mathsf{Cat}$, the following categories are isomorphic (the notation 
of \ref{eq:D-cat}) is used).
\begin{itemize}
\item[{---}] The {\em category of modules} of the polyad $(D^1,d)$ in
  \cite[Section 2.2]{Bru}. Recall that an object consists of objects
  $\{q_x\}$ in $C_x$ for all $x\in D^0$, together with morphisms
  $\{\xymatrix{d(f)q_{s(f)} \ar[r]|-{\,\, \varrho_f\,\,}& q_{t(f)}}\}$ in
  $C_{t(f)}$ for all $f\in D^1$, such that the following diagrams commute for
  all $x\in D^0$ and all $(f,g)\in D^1\circ D^1$. 
$$
\xymatrix@C=35pt{
(d(f)\circ d(g))q_{s(g)} \ar[r]^-{d(f)\varrho_g} \ar[d]_-{(\mu_{f,g})_{q_{s(g)}}} &
d(f)q_{s(f)} \ar[d]^-{\varrho_f} \\
d(f.g)q_{s(g)} \ar[r]_-{\varrho_{f.g}} &
q_{t(f)}}\qquad
\xymatrix{
q_x\ar[d]_-{(\eta_x)_{q_x}} \ar@{=}[rd] \\
d(e_x)q_x \ar[r]_-{\varrho_{e_x}} &
q_x}
$$
A morphism $(q,\varrho) \to (q',\varrho')$ consists of morphisms
$\{\xymatrix{q_x\ar[r]|-{\, \chi_x\,}& q'_x}\}$ in $C_x$ for all $x\in D^0$ such
that the following diagram commutes for all 
$g\in D^1$.
$$
\xymatrix@C=35pt{
d(g)q_{s(g)} \ar[r]^-{d(g)\chi_{s(g)}} \ar[d]_-{\varrho_g} &
d(g)q'_{s(g)} \ar[d]^-{\varrho'_g} \\
q_{t(g)} \ar[r]_-{\chi_{t(g)}} &
q'_{t(g)}}
$$
\item[{---}]
The full subcategory of the Eilenberg-Moore category of the induced monad
$\mathsf{Span}\vert \mathsf{Cat}((1,\mathbf{1}),$ $(D^1,d))$  
on $\mathsf{Span}\vert \mathsf{Cat}((1,\mathbf{1}),(D^0,C))$ whose 
objects are precisely those Eilenberg--Moore algebras 
$((Q,q),(r,\varrho))$ whose underlying span $Q$ is 
$\xymatrix@C=14pt{
D^0 & \ar@{=}[l] D^0 \ar[r]|(.5){\,!\,} & 1}$. 
\end{itemize}

For any monad $(D^1,d)$ on any 0-cell $(D^0,C)$ in $\mathsf{Span}\vert  
\mathsf{Cat}$, also the following categories are isomorphic (where 
the notation of (\ref{eq:D-cat}) is used).
\begin{itemize}
\item[{---}] The {\em category of representations} of the polyad $(D^1,d)$ in
  \cite[Section 2.3]{Bru}. Recall that an object consists of objects $\{W_k\}$
  of $C_{t(k)}$ for all $k\in D^1$ together with morphisms 
$\{\xymatrix@C=30pt{
d(g)W_k\ar[r]|(.52){\,\, \varrho_{g,k}  \,\,} & W_{g.k}}\}$ for $(g,k)\in
D^1\circ D^1$, rendering commutative the following diagrams for all
$f,g,k\in D^1\circ D^1\circ D^1$.
$$
\xymatrix@C=35pt{
(d(f)\circ d(g))W_k \ar[r]^-{d(f)\varrho_{g,k}} \ar[d]_-{(\mu_{f,g})_{W_k}} &
d(f) W_{g.k} \ar[d]^-{\varrho_{f,g.k}} \\
d(f.g) W_k \ar[r]_-{\varrho_{f.g,k}} &
W_{f.g.k}} \qquad
\xymatrix@C=35pt{
W_k \ar[d]_-{(\eta_{t(k)})_{W_k}} \ar@{=}[rd] \\
d(e_{t(k)})W_k \ar[r]_-{\varrho_{e_{t(k)},k}} &
W_k}
$$
A morphism $(W,\varrho) \to (W',\varrho')$ consists of morphisms 
$\{\xymatrix{
W_k \ar[r]|(.45){\, \varphi_k \,} & W'_k}\}$
such that the following diagram commutes for all $(g,k)\in D^1\circ 
D^1$. 
$$
\xymatrix@C=35pt{
d(g)W_k \ar[r]^-{d(g)\varphi_k} \ar[d]_-{\varrho_{g,k}} &
d(g)W'_k \ar[d]^-{\varrho'_{g,k}} \\
W_{g.k} \ar[r]_-{\varphi_{g.k}} &
W'_{g.k}}
$$
\item[{---}]
The following non-full subcategory of the Eilenberg-Moore category of the
monad $\mathsf{Span}\vert \mathsf{Cat}((1,\mathbf{1}),(D^1,d))$ on
$\mathsf{Span}\vert \mathsf{Cat}((1,\mathbf{1}),(D^0,C))$.
The objects are precisely those Eilenberg--Moore algebras $((Q,q),(r,\varrho))$
whose underlying span $Q$ is
$\xymatrix@C=14pt{
D^0 & \ar[l]|(.45){\, t\, } D^1 \ar[r]|(.45){\,!\,} & 1}$
and whose map $r:D^1\wedge D^1 \to D^1$ is the composition in the 
category $D^1$.
The morphisms are those morphisms of Eilenberg--Moore algebras $(f,\varphi)$
in which $f:D^1\to D^1$ is the identity map. 
\end{itemize}

\subsection{Opmonoidal monads in $\mathsf{Span}\vert \mathsf{Cat}$ 
versus opmonoidal polyads} \label{sec:Span|Cat-opmonoidal-monad}
Combining the descriptions in Sections \ref{sec:monad} and 
\ref{sec:Span|OpMon}, 
we obtain coincidence of the following notions. 
\begin{itemize}
\item[{---}] {\em Opmonoidal polyad} in \cite[Paragraph 2.5]{Bru}. 
That is, a pair consisting of a category and -- regarding this category as a 
bicategory with only identity 2-cells -- a lax functor from it to 
$\mathsf{OpMon}$. 
\item[{---}] Monad in $\mathsf{Span}\vert \mathsf{OpMon}$.
\item[{---}] Opmonoidal monad in $\mathsf{Span}\vert 
\mathsf{Cat}$ living on a monoidale of the form in (\ref{eq:opmap-monoidale}). 
\end{itemize}

From the isomorphism in Section \ref{sec:Mnd(Span|OpMon)}, for any given
category (\ref{eq:D-cat}) we have isomorphism of the following bicategories. 
\begin{itemize} 
\item[{---}] The bicategory of opmonoidal polyads over the category
  (\ref{eq:D-cat}) in \cite[Section 3]{Bru} (see the top of its page 
  18). That is, the bicategory of lax functors from 
  (\ref{eq:D-cat}) to $\mathsf{OpMon}$, lax natural
  transformations and modifications.
\item[{---}]
The following locally full sub-bicategory in the bicategory of monads 
(in the sense of \cite{Str:FTMI}) in $\mathsf{Span}\vert \mathsf{OpMon}$. 
The 0-cells are those monads which live on 0-cells $D^{0} \to 
\mathsf{OpMon}^0$ (for the given object set $D^{0}$), whose 1-cells 
are of the form  
$(\xymatrix@C=16pt{
D^0 & D^1 \ar[l]|(.45){\, t \,} \ar[r]|(.45){\, s \,} & D^0}, 
\mathbf{d}:D^1\to \mathsf{OpMon}^1)$ 
(in terms of the given data $s,t$), and whose multiplication and unit 
2-cells have the respective forms $(\cdot,\mu)$ and $(e,\eta)$ (with 
the given maps $\cdot$ and $e$). The 1-cells are those monad 
morphisms $((H,\mathbf{h}),(f,\varphi))$ whose underlying span $H$ is 
the trivial span $D^0=D^0=D^0$ and whose map $f$ is the canonical 
isomorphism $D^0 \circ D^1 \cong D^1 \cong D^1 \circ D^0$. The 
2-cells are all possible monad transformations $(g,\gamma)$ ($g$ in 
them is necessarily the identity map $D^0\to D^0$). 
\item[{---}] The following locally full sub-bicategory in the 
bicategory of monads in the bicategory 
$\mathsf{OpMon}(\mathsf{Span}\vert \mathsf{Cat}$).   
The 0-cells are those monads which live on monoidales 
with object part 
$\xymatrix@C=18pt{
D^0 \ar[r]|(.4){\,C\,} & \mathsf{Cat}^0}$ (for the given object set 
$D^0$) and with multiplication and unit of the form 
$$
(\xymatrix@C=16pt{
D^0 =D^0 \ar[r]^-\Delta & D^0 \times D^0},
\xymatrix@C=16pt{
D^0 \ar[r]^-\otimes & \mathsf{Cat}^1})
\quad \textrm{and} \quad
(\xymatrix@C=16pt{
D^0=D^0 \ar[r]^-{!} & 1},
\xymatrix@C=16pt{
D^0 \ar[r]^-K & \mathsf{Cat}^1}),
$$
whose 1-cells are of the form 
$(\xymatrix@C=16pt{
D^0 & D^1 \ar[l]|(.45){\, t \,} \ar[r]|(.45){\, s \,} & D^0}, 
d:D^1 \to \mathsf{Cat}^1)$ 
(in terms of the given data $s,t$), and whose multiplication and unit 
2-cells have the respective forms $(\cdot,\mu)$ and $(e,\eta)$ (with 
the given maps $\cdot$ and $e$). The 1-cells are those monad 
morphisms $((H,h),(f,\varphi))$ whose underlying span $H$ is the trivial 
span $D^0=D^0=D^0$ and whose map $f$ is the canonical isomorphism 
$D^0 \circ D^1 \cong D^1 \cong D^1 \circ D^0$. The 2-cells are all 
possible monad transformations $(g,\gamma)$ ($g$ in them is 
necessarily the identity map $D^0\to D^0$).  
\end{itemize}

\subsection{Hopf monads in $\mathsf{Span}\vert \mathsf{Cat}$ versus Hopf
  polyads} 
Our next task is to compute the fusion 2-cells as in \cite{CLS} for 
the opmonoidal monads in $\mathsf{Span}\vert \mathsf{Cat}$ of Section 
\ref{sec:Span|Cat-opmonoidal-monad}.  
The left fusion 2-cell consists of the map of spans
\begin{eqnarray*}
&&{\color{white}.}\hspace{-.7cm}
(\xymatrix@C=8pt{
D^0 & \ar[l] \{(p,q)\in D^1 \times D^1 \vert s(p)=t(q)\} \ar[r] & D^0\times D^0,}
t(p) \mapsfrom (p,q) \mapsto (s(q),s(p))) \to \\
&&{\color{white}.}\hspace{-.7cm}
(\xymatrix@C=8pt{
D^0 & \ar[l] \{(p,q)\in D^1 \times D^1 \vert t(p)=t(q)\} \ar[r] & D^0\times D^0,}
t(p) \mapsfrom (p,q) \mapsto (s(p),s(q)))
\end{eqnarray*}
sending $(p,q)$ to $(p.q,p)$; and the set of natural 
transformations 
\begin{equation} \label{eq:Span|Cat-Hopf}
\xymatrix@C=20pt{
d(p)(d(q)(-) 
\raisebox{-5pt}{$
\stackrel{\displaystyle \otimes}{{}_{s(p)}}
$}
(-)) 
\ar[r]^-{\raisebox{10pt}{$_{d^2_p}$}} &
(d(p)\circ d(q))(-) 
\raisebox{-5pt}{$
\stackrel{\displaystyle \otimes}{{}_{t(p)}}
$}
d(p)(-) 
\ar[rr]^-{\raisebox{10pt}{$_{\mu_{p,q} \!
\raisebox{-3pt}{$
\stackrel{\otimes}{{}_{{}_{t(p)}}}
$}
\! 1}$}} &&
d(p.q)(-) 
\raisebox{-5pt}{$
\stackrel{\displaystyle \otimes}{{}_{t(p)}}
$}
d(p)(-)} 
\end{equation}
between functors $C_{s(q)}\times C_{s(p)} \to C_{t(p)}$,
labelled by $(p,q) \in D^1 \wedge D^1$ (a label $x\in D^0$ on 
$\otimes$ refers to the category $C_x$ in which it serves as the 
monoidal product). This coincides with the left fusion operator of 
\cite[Definition 2.15]{Bru}.  

Clearly, this left fusion 2-cell above is invertible in $\mathsf{Span}\vert 
\mathsf{Cat}$ if and only if the underlying category (\ref{eq:D-cat}) is a
groupoid and each natural transformation in the set 
(\ref{eq:Span|Cat-Hopf}) is invertible.
So we obtained the coincidence of the following notions.
\begin{itemize}
\item[{---}] {\em Left Hopf polyad} in the sense of \cite[Definition
2.17]{Bru} whose underlying category is a groupoid. That is, an opmonoidal
  polyad whose underlying category is a groupoid and for which each of the
  natural transformations (\ref{eq:Span|Cat-Hopf}) is invertible.
\item[{---}] A Hopf monad in $\mathsf{Span}\vert \mathsf{Cat}$ living on a
monoidale of the form in (\ref{eq:opmap-monoidale}). 
\end{itemize}
The case of the right fusion 2-cell is symmetric.

\subsection{The induced Hopf monad in $\mathsf{Cat}$} 
Since the monoidal product in $\mathsf{Cat}$ is Cartesian, any 0-cell (that
is, any category) is a comonoidale in a unique way. Hence the construction in
Section \ref{sec:induced_monoidale} yields an induced comonoidale $(Y,C)$ in  
$\mathsf{Span}\vert \mathsf{Cat}$ for any set of categories 
$\{C_y\vert y\in Y\}$. 

On the other hand, as described in Section 
\ref{sec:Span|OpMon}, any set of {\em monoidal} categories 
$\{(M_x,\otimes_x,$ $K_x)\vert x\in X\}$ induces a monoidale $(X,M)$ in  
$\mathsf{Span}\vert \mathsf{Cat}$. So there is a monoidal 
category $\mathsf{Span}\vert \mathsf{Cat}((Y,C),(X,M))$ as in Section 
\ref{sec:conv_mod}.

Let $(D^1,d)$ be an opmonoidal polyad on $(D^0,M)$; that is, an 
opmonoidal monad in $\mathsf{Span}\vert \mathsf{Cat}$. It induces an 
opmonoidal monad in $\mathsf{Cat}$ on the category 
$\mathsf{Span}\vert \mathsf{Cat}((Y,C),$ $(D^0,M))$, see again Section 
\ref{sec:conv_mod}. One can define its {\em Hopf modules} as in \cite{BV} and
\cite[Section 6.5]{BLV}. Criteria for the equivalence between the category of
these Hopf modules and $\mathsf{Span}\vert \mathsf{Cat}((Y,C),(D^0,M))$ were
obtained in \cite[Theorem 6.11]{BLV}; known as the {\em fundamental theorem of 
Hopf modules}. 

The inclusion of the category of representations of a polyad 
into the Eilenberg-Moore category of the induced monad in Section 
\ref{sec:Cat_induced_monad} lifts to an inclusion of
the category of {\em Hopf representations} in \cite[Section 6.2]{Bru} 
into the above category of Hopf modules in the sense of \cite{BLV},
for $(Y,C)=(1,\mathbf{1})$. Hence if the fundamental theorem 
of Hopf modules in \cite{BLV} holds, then the equivalence therein 
induces an equivalence between this subcategory in \cite[Section 
6.2]{Bru} and a suitable subcategory of $\mathsf{Span}\vert 
\mathsf{Cat}((1,\mathbf{1}),(D^0,M))$. This gives an alternative proof 
of \cite[Theorem 6.3]{Bru}. 

On the other hand, the category of {\em Hopf modules} in \cite[Section 
6.1]{Bru} does not seem to be a subcategory of the above category of Hopf
modules in the sense of \cite{BLV}, for $(Y,C)=(1,\mathbf{1})$; and
\cite[Theorem 6.1]{Bru} seems to be of different nature.

\section{Hopf group monoids and Hopf categories as Hopf monads on 
naturally Frobenius opmap monoidales} 
For an arbitrary object $X$ in any bicategory $\mathcal M$, a {\em monad on}
$X$ is exactly the same thing as a {\em monoid in} the monoidal endohom
category $\mathcal M(X,X)$ --- though one of these equivalent descriptions may
turn out to be more convenient in one or another situation.

If $X$ is an {\em opmap monoidale} (that is, a monoidale or pseudo-monoid whose 
multiplication and unit 1-cells possess left adjoints) in a {\em monoidal
bicategory} $\mathcal M$, then the endohom category $\mathcal M(X,X)$ possesses
the richer structure of a so-called {\em duoidal category}; see \cite{St}.

A {\em duoidal} (or {\em 2-monoidal} in the terminology of 
\cite{AgMa}) category is a category with two monoidal structures 
$(\circ, I)$ and $(\bullet, J)$ which are compatible in the sense that 
the functors $\circ$ and $I$, as well as their associativity and unitality 
natural isomorphisms are opmonoidal for the $\bullet$-product. 
Equivalently, the functors $\bullet$ and $J$, as well as their associativity
and unitality natural isomorphisms are monoidal for the $\circ$-product. In
technical terms it means the existence of four natural transformations (the
binary and nullary parts of two (op)monoidal functors) subject to a number of
conditions spelled out e.g. in \cite{AgMa}. 

For an opmap monoidale $X$ in a monoidal bicategory $\mathcal M$, the first
monoidal product $\circ$ on $\mathcal M(X,X)$ comes from the horizontal
composition $\circ$ in $\mathcal M$. Since $X$ possesses both structures of a
monoidale and a comonoidale (the latter one with the comultiplication and
the counit provided by the adjoints of the multiplication and the unit),
$\mathcal M(X,X)$ has a second monoidal product $\bullet$ of the convolution
type, see Section \ref{sec:conv_mod}. Thanks to the (adjunction) relation
between the monoidale and the comonoidale $X$, these monoidal structures
$\circ$ and $\bullet$ render $\mathcal M(X,X)$ with the structure of duoidal
category. 

This observation turns out to be very useful: the coincidence 
of a {\em monad on} $X$ and a {\em monoid in} $(\mathcal M(X,X),\circ)$ is 
supplemented with the coincidence of an {\em opmonoidal endo 1-cell 
on} $X$ and a {\em comonoid in} $(\mathcal M(X,X),\bullet)$; see 
\cite[Section 3.3]{BoLa}. Combining these correspondences, an {\em 
opmonoidal monad on} an opmap monoidale $X$ in a monoidal bicategory 
$\mathcal M$ turns out to be exactly the same thing as a {\em 
bimonoid in} the duoidal endohom category $\mathcal M(X,X)$ (in the 
sense of \cite[Definition 6.25]{AgMa}), see again \cite{St} or a 
review  in \cite[Section 3.3]{BoLa}. 

Although these are {\em mathematically equivalent} points of view, 
one of them may turn out to be {\em more convenient} in one or 
another situation. Recall for example, that no sensible notion of 
{\em antipode} for Hopf monads {\em on} arbitrary monoidales of 
monoidal bicategories is known. It is one of the key observations in 
\cite{BoLa}, however, that for a Hopf monad living on a {\em 
naturally Frobenius opmap monoidale}, it can be given a natural 
meaning. In this situation, the antipode axioms are formulated most 
easily {\em in} the duoidal endohom category, see \cite[Theorem 
7.2]{BoLa}.

Since in this section we shall study Hopf-like structures --- {\em Hopf group 
monoids} and {\em Hopf categories} --- defined in terms of antipode 
morphisms, we are to apply this language.  

A braided monoidal small category $(V,\otimes, K,c)$ can be regarded 
as a monoidal bicategory with a single object, in this section 
we will work with that.

\subsection{The bicategory $\mathsf{OpMon}(V)$ for a braided monoidal 
category $V$} \label{sec:OpMon(V)}

An object of $\mathsf{OpMon}(V)$ --- that is, a monoidale in $V$ --- 
consists of two objects $M$ and $U$ of $V$ (the multiplication and 
the unit) and three coherence isomorphisms $\alpha:M\otimes M \to 
M\otimes M$, $\lambda:M\otimes U \to K$ and $\varrho:M\otimes U \to 
K$ subject to the appropriate pentagon and triangle conditions. 

Here we are not interested in arbitrary monoidales in $V$. The one 
which plays a relevant role is the {\em trivial} one which has both 
the multiplication and the unit equal to the monoidal unit $K$ and 
all coherence isomorphisms built up from the coherence isomorphisms 
of $V$.

A 1-cell of $\mathsf{OpMon}(V)$ --- that is, an opmonoidal 1-cell 
in $V$ ---  is an object $A$ of $V$ equipped with morphisms $a^2: 
A\otimes M \to M'\otimes A \otimes A$ and $a^0: A\otimes U \to U'$ 
subject to appropriate coassociativity and counitality conditions.

The endo 1-cells of the trivial monoidale are then the same as the 
comonoids $(A,a^2,a^0)$ in $V$.

A 2-cell of $\mathsf{OpMon}(V)$ --- that is, an opmonoidal 2-cell 
in $V$ --- is a morphism $A\to A'$ in $V$ which is appropriately 
compatible with the opmonoidal structures $(a^2,a^0)$ and $(a^{\prime 
2},a^{\prime 0})$.

Between endo 1-cells of the trivial monoidale, the 2-cells are then 
the same as the comonoid morphisms $(A,a^2,a^0)\to (A',a^{\prime 
2},a^{\prime 0})$.

So for any braided monoidal category $V$, we obtain isomorphism of the
following monoidal categories. 
\begin{itemize}
\item[{---}] The endohom category of the trivial monoidale in 
$\mathsf{OpMon}(V)$. 
\item[{---}] The category $\mathsf{Cmd}(V)$ of comonoids in $V$.   
\end{itemize}    

\subsection{Sets as naturally Frobenius opmap monoidales in $\mathsf{Span}\vert
  V$} \label{sec:Span|V-monoidale}
Since there is only one 0-cell of the bicategory $V$, the 0-cells of 
$\mathsf{Span}\vert V$ are simply sets. Moreover, the only 0-cell of 
the bicategory $V$ is the monoidal unit, hence it is a trivial 
monoidale, so in particular a naturally Frobenius opmap monoidale. 
Thus for any set $X$ the construction in Section \ref{sec:induced_monoidale}
yields a naturally Frobenius opmap monoidale in $\mathsf{Span}\vert V$ with
multiplication and unit 1-cells consisting of the respective spans 
$$
\xymatrix{
X & X \ar@{=}[l] \ar[r]^-\Delta  & X\times X}
\quad
\textrm{and}
\quad
\xymatrix{
X & X \ar@{=}[l] \ar[r]^-{!}  & 1}
$$
and in both cases the constant map sending each element of $X$ to 
the monoidal unit $K$ of $V$; and trivial (i.e. built up from coherence
isomorphisms of $V$) associativity and unitality coherence 2-cells.

\subsection{The bicategory $\mathsf{Span}\vert \mathsf{OpMon}(V)$} 
\label{sec:Span|OpMon(V)}
The isomorphism of Section \ref{sec:Mnd(Span|OpMon)} takes an object of 
$\mathsf{OpMon}(\mathsf{Span}\vert V)$ of the form in Section
\ref{sec:Span|V-monoidale} to the object of $\mathsf{Span}\vert
\mathsf{OpMon}(V)$ which consists of the set $X$ and the constant map sending
each element of $X$ to the trivial monoidale in $V$ (see Section
\ref{sec:OpMon(V)}). For brevity we will denote simply by $X$ also this object
of $\mathsf{Span}\vert \mathsf{OpMon}(V)$. We are interested in the 2-full
sub-bicategory of $\mathsf{Span}\vert \mathsf{OpMon}(V)$ defined by these
objects.  

For any sets $X$ and $Y,$ an object of $\mathsf{Span}\vert 
\mathsf{OpMon}(V)(X,Y)$ consists of a span
$\xymatrix@C=7pt{
Y & \ar[l] A \ar[r] & X}$
and a map from $A$ to the object set of the endohom category of the 
trivial monoidale in $\mathsf{OpMon}(V)$. That is, in view of the 
isomorphism of Section \ref{sec:OpMon(V)}, a map $a$ from $A$ to the 
set of comonoids in $V$.
 
A morphism in $\mathsf{Span}\vert \mathsf{OpMon}(V)(X,Y)$ consists of 
a map of spans $f:A\to A'$ and morphisms $a(p)\to a'f(p)$ in the 
endohom category of the trivial monoidale in $\mathsf{OpMon}(V)$, 
for all $p\in A$. That is, in view of Section \ref{sec:OpMon(V)}, a set of
comonoid morphisms $\{a(p)\to a'f(p)\ \vert\ p\in A\}$ in $V$. 

This leads to an isomorphism between the following categories, for 
any sets $X,Y$ and any braided monoidal category $V$.
\begin{itemize}
\item[{---}] $\mathsf{OpMon}(\mathsf{Span}\vert V)(X,Y)$.
\item[{---}] $\mathsf{Span}\vert \mathsf{OpMon}(V)(X,Y)$.    
\item[{---}] $\mathsf{Span}\vert\mathsf{Cmd}(V)(X,Y)$.
\end{itemize}    

\subsection{The duoidal endohom categories} 
\label{sec:Span|V-endohom}

The structure of an opmap monoidale that we constructed in Section
\ref{sec:Span|V-monoidale} on any set $X$, induces a duoidal structure on the
endohom category $\mathsf{Span}\vert V(X,X)$ which we describe next. 
It is obtained by a straightforward application of the general construction in
\cite{St}, see also \cite[Section 3.3]{BoLa}.   

The objects of $\mathsf{Span}\vert V(X,X)$ are pairs consisting of an 
$X$-span $A$ and a map $a$ from the set $A$ to the set of objects in $V$. 
The morphisms $(A,a)\to (A',a')$ are pairs consisting of a map of 
$X$-spans $f:A\to A'$ and a set $\{\varphi_h:a(h) \to a'f(h) \vert h\in A\}$
of morphisms in $V$.   

The first monoidal product $\circ$ on $\mathsf{Span}\vert V(X,X)$ comes from
the horizontal composition in $\mathsf{Span}\vert V$; thus in fact from the
monoidal product in $V$: the product of any two morphisms $(g,\gamma): (B,b)
\to (B',b')$ and $(f,\varphi): (A,a) \to (A',a')$
is 
$$
(g\wedge f:B\wedge A \to B' \wedge A',
\{ \gamma_d \otimes \varphi_h:b(d) \otimes f(h) \to b'g(d) \otimes 
a'f(h) \vert (d,h) \in B\wedge A \}).
$$
The monoidal unit $I$ is the identity 1-cell of $X$: it consists of the 
trivial $X$-span and the map sending each element of $X$ to the 
monoidal unit $K$ of $V$. 

For any (possibly different) opmap monoidales $X$ and $Y$ of the kind 
discussed in Section \ref{sec:Span|V-monoidale}, the hom category 
$\mathsf{Span} \vert V(X,Y)$ admits a monoidal product $\bullet$ 
which is of the convolution type, see Section \ref{sec:conv_mod}.
Now the product of 2-cells 
$(g,\gamma): (B,b)\Rightarrow (B',b')$ and $(f,\varphi): (A,a) 
\Rightarrow (A',a')$ between 1-cells $X \to Y$ is the pair
consisting of the map of spans in (\ref{eq:prod_span}) and the set 
$\{ \gamma_d \otimes \varphi_h: b(d)\otimes a(h) \to b'g(d) \otimes 
a'f(h) \vert (d,h) \in B \bullet A \}$ of morphisms in $V$.
The monoidal unit $J$ consists of the complete span
$\xymatrix@C=12pt{
Y & Y \times X \ar[l] \ar[r] & X}$ 
and the map sending each element of $Y\times X$ to the monoidal unit 
$K$ of $V$.

The above monoidal structures combine into a duoidal structure on
$\mathsf{Span}\vert V(X,X)$.
The four structure morphisms take the following forms. The first one 
is a morphism
$((A,a) \bullet (B,b)) \circ ((H,h) \bullet (D,d)) \to
((A,a)\circ (H,h)) \bullet ((B,b)\circ (D,d))$
which is natural in each object $(A,a),(B,b),(H,h),(D,d)$. It 
consists of the map of spans
$$
\xymatrix@C=1pt{
& \raisebox{8pt}{$
\makebox[.7\textwidth][l]{\hspace{-2.2cm}
$\{(p,q,v,w) \in A\times B \times H\times D \vert
l(p) = l(q),\ r(p)=r(q)=l(v)=l(w),\ r(v)=r(w)\}$}$}
\ar[ld]_-{(p,q,v,w) \mapsto l(p)\quad}
\ar[rd]^-{\quad (p,q,v,w) \mapsto r(v)}
\ar@{>->}[dd]^-{(p,q,v,w) \mapsto (p,v,q,w)}\\
X && X \\
& \makebox[.7\textwidth][l]{\hspace{-2cm}
$\{(p,v,q,w) \in A\times H \times B\times D \vert
l(p) = l(q),\ r(p)=l(v),\ r(q)=l(w),\ r(v)=r(w)\}$}
\ar[lu]^-{(p,v,q,w) \mapsto l(p)\quad}
\ar[ru]_-{\quad (p,v,q,w) \mapsto r(v)}}
$$
and the set
$$
\{1\otimes c \otimes 1: a(p) \otimes b(q) \otimes h(v) \otimes d(w) \to
a(p)  \otimes h(v) \otimes b(q) \otimes d(w) \}
$$
of morphisms in $V$, labelled by the elements $(p,q,v,w)\in (A\bullet B)\wedge
(H \bullet D)$. 

Next we need a morphism $J\circ J \to J$; it consists of
the map of spans
$$
\xymatrix@C=50pt{
& X \times X \times X 
\ar[ld]_-{(p,q,v) \mapsto p} \ar[rd]^-{(p,q,v) \mapsto v} 
\ar[dd]^-{(p,q,v) \mapsto (p,v)} \\
X && X \\
& X\times X 
\ar[lu]^-{(p,q) \mapsto p} \ar[ru]_-{(p,q) \mapsto q}
}
$$
and the map sending each element of $X \times X \times X$ to the identity
morphism of the monoidal unit $K$ of $V$.

Then we need a morphism $I \to I \bullet I=I$; it is the identity morphism.

Finally we need a morphism $I\to J$. It is given by the diagonal map 
$\Delta:X \to X\times X$ from the trivial to the complete span and 
the map sending each element of $X$ to the identity morphism of the 
monoidal unit $K$ of $V$. 

\subsection{The Zunino category}
There is a particular duoidal category $\mathsf{Span}\vert V(1,1)$ of 
the above form in Section \ref{sec:Span|V-endohom} for the singleton 
set $1$. Here both monoidal products $\circ$ and $\bullet$ turn out 
to be equal, and sending any pair of 2-cells $(g,\gamma): (B,b) 
\Rightarrow (B',b')$ and $(f,\varphi): (A,a) \Rightarrow (A',a')$ 
between 1-cells $1\to 1$ to 
$$
(g\times f:B\times A \to B' \times A',
\{ \gamma_d \otimes \varphi_h:b(d) \otimes f(h) \to b'g(d) \otimes 
a'f(h) \vert (d,h) \in B\times A \}).
$$
This amounts to saying that the duoidal category $\mathsf{Span}\vert 
V(1,1)$ coincides with the braided monoidal {\em Zunino category}; for its
explicit description (in the case when $V$ is the symmetric monoidal category
of modules over a commutative ring) see \cite[Section 2.2]{CL}.

\subsection{Hopf group monoids}
For an ordinary monoid $G$ (that is, a monoid in the Cartesian monoidal
category of sets), a {\em $G$-algebra} was defined in \cite[Definition
1.6]{CL} as a monoidal functor from $G$ --- regarded as a discrete category
with object set $G$ and monoidal structure coming from the multiplication
$\cdot$ and unit $e$ of $G$ --- to the monoidal category of vector spaces
(over a given field). Following this idea, we define a {\em $G$-monoid} in any
monoidal category $V$ as a monoidal functor from $G$ to $V$. This is the same
as a lax functor from the 1-object category $G$ (regarded as a bicategory with
only identity 2-cells) to $V$ (regarded as a bicategory with a single
0-cell). Hence from Section \ref{sec:monad}, and from the correspondence
between monads on some object and monoids in its composition-monoidal endohom
category, we obtain the coincidence of the following notions for any monoidal
category $V$. 
\begin{itemize}
\item[{---}] A pair consisting of an ordinary monoid $G$ and a {\em
  $G$-monoid} in $V$. 
\item[{---}] A monad in $\mathsf{Span}\vert V$ on the singleton set $1$.
\item[{---}] A monoid in the Zunino category $\mathsf{Span}\vert 
V(1,1)$.
\end{itemize}

Combining the isomorphism of Section \ref{sec:Span|OpMon(V)}, and the
correspondence of opmonoidal 1-cells on some opmap monoidale and
comonoids in its convolution-monoidal endohom category, the following
categories are isomorphic for any braided monoidal category $V$. 
\begin{itemize}
\item[{---}] The endohom category of the singleton set $1$ in 
$\mathsf{Span}\vert \mathsf{Cmd}(V)$.
\item[{---}] The endohom category of the singleton set $1$ --
regarded as an opmap monoidale in Section \ref{sec:Span|V-monoidale} -- in 
$\mathsf{OpMon}(\mathsf{Span}\vert V)$.
\item[{---}] The category of comonoids in the Zunino category 
$\mathsf{Span}\vert V(1,1)$.
\end{itemize}    

For any monoid $G$, a {\em semi Hopf $G$-algebra} was defined in 
\cite[Definition 1.7]{CL} as a $G$-monoid (in the above sense) in the monoidal
category of coalgebras (over a given field). Following this idea, we define a
{\em semi Hopf $G$-monoid} in any braided monoidal category $V$ as a
$G$-monoid in $\mathsf{Cmd}(V)$. Hence combining the isomorphism above, and
the correspondence between opmonoidal monads on some opmap monoidale and
bimonoids in its duoidal endohom category, we obtain the coincidence of the
following notions for any monoid monoidal category $V$.   
\begin{itemize}
\item[{---}] A pair consisting of a monoid $G$ and a semi Hopf 
$G$-monoid in $V$. 
\item[{---}] A monad in $\mathsf{Span}\vert \mathsf{Cmd}(V)$ on the 
singleton set $1$. 
\item[{---}] An opmonoidal monad in $\mathsf{Span}\vert V$ on 
the monoidale $1$.
\item[{---}] A bimonoid in the Zunino category $\mathsf{Span}\vert 
V(1,1)$.
\end{itemize}    

For a group $G$, a semi Hopf $G$-algebra --- that is, a  monoidal 
functor from the discrete category on the object set $G$ to the monoidal
category of coalgebras, sending $p\in G$ to a coalgebra
$(g(p),\delta_p,\varepsilon_p)$; with binary part of the monoidal structure
denoted by   
$\{\xymatrix@C=32pt{ 
g(p)\otimes g(q) \ar[r]|-{\,\,\mu_{p,q}\,\,} &
g(p.q)}\}_{p,q\in G}$
and nullary part denoted by 
$\xymatrix@C=18pt{
K \ar[r]|-{\,\eta\,} & g(e)}$ ---
was termed a {\em Hopf $G$-algebra} in \cite[Definition 1.8]{CL} if 
equipped with linear maps (the so-called {\em antipode}) 
$\{\xymatrix{ 
g(p) \ar[r]|-{\,\sigma_p\,} &
g(p^{-1})}\}_{p\in G}$ rendering commutative the following diagram for all
$p\in G$.
$$
\xymatrix@R=10pt{
g(p) \ar[r]^-{\delta_p} \ar[dd]_-{\delta_p} \ar[rd]_-{\varepsilon_p}&
g(p) \otimes g(p) \ar[r]^-{\sigma_p \otimes 1}  &
g(p^{-1})\otimes g(p) \ar[dd]^-{\mu_{p^{-1},p}} \\
& K \ar[rd]^-\eta \\
g(p) \otimes g(p) \ar[r]_-{1\otimes \sigma_p} &
g(p) \otimes g(p^{-1}) \ar[r]_-{\mu_{p,p^{-1}}} &
g(e).}
$$ 
By this motivation we define a {\em Hopf $G$-monoid} in any braided 
monoidal category $V$ as a monoidal functor $((g,\delta,\varepsilon),
\mu,\eta)$ from the discrete category on the object set $G$ to
$\mathsf{Cmd}(V)$ together with morphisms  
$\{\xymatrix{ 
g(p) \ar[r]|-{\,\sigma_p\,} &
g(p^{-1})}\}_{p\in G}$ 
in $V$ rendering commutative the same diagram.

Note that this diagram encodes precisely the antipode axioms of \cite[Theorem
7.2]{BoLa} for the bimonoid $g$ in the duoidal Zunino category
$\mathsf{Span}\vert V(1,1)$; which are in turn the same as the usual antipode
axioms for the bimonoid $g$ in the braided monoidal Zunino category
$\mathsf{Span}\vert V(1,1)$. 
Thus since the singleton set is regarded as a naturally Frobenius opmap 
monoidale in $\mathsf{Span}\vert V$ (in the way described in Section 
\ref{sec:Span|V-monoidale}), from \cite[Theorem 7.2]{BoLa} we deduce 
the coincidence of the following notions for any braided monoidal 
category $V$.
\begin{itemize}
\item[{---}] A pair consisting of a group $G$ and a Hopf $G$-monoid in $V$.
\item[{---}] A Hopf monoid in the Zunino category $\mathsf{Span}\vert 
V(1,1)$.
\item[{---}] A Hopf monad in $\mathsf{Span}\vert V$ on the monoidale $1$. 
\end{itemize} 
 
\subsection{Monads in $\mathsf{Span}\vert V$ versus categories 
enriched in $V$} \label{sec:Span|V-monad} 
We turn to the interpretation of $V$-enriched categories in \cite[Section
2]{BCV} as monads in $\mathsf{Span}\vert V$, matrices of comonoids in $V$ as
in \cite[Section 3]{BCV} as opmonoidal 1-cells in $\mathsf{Span}\vert V$, 
categories enriched in the category of comonoids in $V$ as in
\cite[Proposition 3.1]{BCV} as opmonoidal monads in $\mathsf{Span}\vert V$,
and finally the Hopf categories of \cite[Definition 3.3]{BCV} as Hopf monads
in $\mathsf{Span}\vert V$.

Recall that a category enriched in $V$ can be described as a pair consisting
of a set $X$ (it plays the role of the set of objects) and a lax functor from
the indiscrete category on the object set $X$, regarded as a bicategory with
only identity 2-cells, to $V$, regarded as a bicategory with a single object.
An identity-on-objects $V$-enriched functor is precisely a lax natural
transformation whose 1-cell part is trivial. 

On the other hand, between monads on the same object in any bicategory, a
monad morphism (in the sense of \cite{Str:FTMI}) with trivial 1-cell part is
precisely the same thing as a morphism between the corresponding monoids in
the composition-monoidal endohom category.

Using these observations and the fact that the complete span 
$\xymatrix@C=12pt{
X & \ar[l] X\times X \ar[r] & X}$
is terminal in $\mathsf{Span}(X,X)$, from Section \ref{sec:lax_funtor_bicat}
we obtain isomorphism of the following categories, for any braided monoidal
category $V$ and any set $X$. 
\begin{itemize}
\item[{---}] The category whose objects are the $V$-enriched categories 
  with object set $X$,
  and whose morphisms are the identity-on-object $V$-enriched functors. (This
  category is used in \cite{BCV}, see its page 1176.)
\item[{---}] The category whose objects are those monads on $X$ in
  $\mathsf{Span}\vert V$ which live on such 1-cells of
  $\mathsf{Span}\vert V$ whose underlying $X$-span is the complete span
$\xymatrix@C=8pt{
X & \ar[l] X\times X \ar[r] & X}$; 
and whose morphisms are those monad morphisms in
  $\mathsf{Span}\vert V$ (in the sense of \cite{Str:FTMI}) whose 1-cell part
is the identity 1-cell $X\to X$ in $\mathsf{Span}\vert V$.  
\item[{---}] The full subcategory of the category of monoids in 
  $(\mathsf{Span}\vert V(X,X),\circ, I)$ whose objects live on such 1-cells of
  $\mathsf{Span}\vert V$ in which the underlying $X$-span is the complete span
$\xymatrix@C=12pt{
X & \ar[l] X\times X \ar[r] & X}$.  
\end{itemize}

\subsection{Opmonoidal 1- and 2-cells in $\mathsf{Span}\vert V$ versus 
matrices of comonoids, and of comonoid morphisms in $V$} 
\label{sec:Span|V-opmonoidal} 
Again, we are not interested in arbitrary opmonoidal 1- and 2-cells only in
those between opmap monoidales $X$ and $Y$ of the kind discussed in 
Section \ref{sec:Span|V-monoidale}. 

Let us use again the fact that the complete span 
$\xymatrix@C=12pt{
Y & \ar[l] Y\times X \ar[r] & X}$
is terminal in $\mathsf{Span}(X,Y)$.
Then from the isomorphism of Section \ref{sec:Span|OpMon(V)} on the one hand,
and from the correspondence between opmonoidal 1-cells on some opmap monoidale
and comonoids in its convolution-monoidal endohom category on the other hand,
we obtain the following isomorphism of full subcategories, for any braided
monoidal category $V$ and any sets $X,Y$.
\begin{itemize}
\item[{---}] The category whose objects are matrices of comonoids in $V$ with
  columns labelled by the elements of $X$ and rows labelled by the elements
  of $Y$; and whose morphisms are $X$ by $Y$ matrices of comonoid
  morphisms in $V$. 
\item[{---}] The full subcategory of opmonoidal 1-cells $X\to Y$ in
  $\mathsf{Span}\vert V$ and opmonoidal 2-cells between them, for whose
  objects the underlying span is the complete span 
$\xymatrix@C=12pt{
Y & \ar[l] Y\times X \ar[r] & X}$.
\item[{---}] The full subcategory of comonoids in $(\mathsf{Span}\vert
  V(X,Y),\bullet, J)$ for whose objects the underlying span is the complete
  span 
$\xymatrix@C=12pt{ 
Y & \ar[l] Y\times X \ar[r] & X}$.
\end{itemize}

\subsection{Opmonoidal monads in $\mathsf{Span}\vert V$ versus 
categories enriched in $\mathsf{Cmd}(V)$} 
\label{sec:Span|V-bimonad} 
From the isomorphisms of Section \ref{sec:Span|V-monad} and Section
\ref{sec:Span|OpMon(V)} on the one hand, and the correspondence between
opmonoidal monads on an opmap monoidale and the bimonoids in its duoidal
endohom category on the other hand, isomorphism of the following categories
follows, for any set $X$ and any braided monoidal category $V$.
\begin{itemize}
\item[{---}] The category whose objects are the $\mathsf{Cmd}(V)$-enriched
categories with object set $X$; and whose morphisms are the identity-on-object
$\mathsf{Cmd}(V)$-enriched functors. (This category is used in \cite{BCV}, see
its page 1177.)  
\item[{---}] The category in which the objects are those opmonoidal monads in
  $\mathsf{Span}\vert V$ on the opmap monoidale $X$ of Section
  \ref{sec:Span|V-monoidale} in whose 1-cell part $X\to X$ the underlying
  span is the complete span   
$\xymatrix@C=12pt{
X & \ar[l] X\times X \ar[r] & X}$; 
and whose morphisms are those opmonoidal monad morphisms whose 1-cell part is
the identity 1-cell $X\to X$ in $\mathsf{OpMon}(\mathsf{Span}\vert V)$.  
\item[{---}] The full subcategory of the category of bimonoids (in the sense
  of \cite[Definition 6.25]{AgMa}) in the duoidal category
  $\mathsf{Span}\vert V(X,X)$, defined by those objects which live on 1-cells
  $X\to X$ in $\mathsf{Span}\vert V$ with underlying span the complete span
$\xymatrix@C=12pt{
X & \ar[l] X\times X \ar[r] & X}$.
\end{itemize}    

\subsection{The induced opmonoidal monad in $\mathsf{Cat}$}
Regard a $V$-enriched category with object set $X$ as a 
monad in $\mathsf{Span}\vert V$ on the 0-cell $X$ as in Section
\ref{sec:Span|V-monad}. Via horizontal composition it induces a monad in
$\mathsf{Cat}$ on the category $\mathsf{Span}\vert V(Y,X)$ for any set $Y$,
see Section \ref{sec:conv_mod}.

If we start with a category enriched in the category of comonoids in $V$ 
--- that is, as a monad in $\mathsf{Span}\vert V$ it admits an opmonoidal 
structure with respect to the monoidale of Section
\ref{sec:Span|V-monoidale}, see Section \ref{sec:Span|V-bimonad} --- then so
does the induced monad   
in $\mathsf{Cat}$ with respect to the convolution monoidal structure of
$\mathsf{Span}\vert V(Y,X)$, see again Section \ref{sec:conv_mod}. This
implies the monoidality (via the product $\bullet$) of the Eilenberg-Moore 
category of the induced monad. 

Consider a $\mathsf{Cmd}(V)$-enriched category with object set $X$ and hom
objects $(a(x,y),$ $\delta_{x,y},\varepsilon_{x,y})$ for $(x,y)\in X \times
X$. Denote the composition compatibility morphisms by $\mu_{x,y,z}:a(x,y)
\otimes a(y,z)\to a(x,z)$ and denote the unit compatibility morphisms by
$\eta_x:K \to a(x,x)$, for all $x,y,z\in X$. For these data, the following
monoidal categories are isomorphic.  
\begin{itemize}
\item[{---}] The {\em category of modules} in \cite[Definition
  4.1]{BCV}. Recall that its objects are sets $\{v(p,q)\}_{p,q\in X}$ of
  objects in $V$ together with sets of morphisms in $V$
$\{\xymatrix@C=35pt{
a(x,y) \otimes v(y,z) \ar[r]|-{\,\psi_{x,y,z}\,} &
v(x,z)}\}_{x,y,z\in X}$ 
making commutative for all $x,y,z,$ $u\in X$ the following
associativity and unitality diagrams.  
$$
\xymatrix@C=35pt{
a(x,y) \otimes a(y,z) \otimes v(z,u) 
\ar[r]^-{\mu_{x,y,z}\otimes 1} \ar[d]_-{1\otimes \psi_{y,z,u}} &
a(x,z) \otimes v(z,u) \ar[d]^-{\psi_{x,z,u}} \\
a(x,y) \otimes v(y,u) \ar[r]_-{\psi_{x,y,u}} &
v(x,u)}\qquad
\xymatrix{
v(x,y) \ar[r]^-{\eta_x \otimes 1} \ar@{=}[rd] &
a(x,x) \otimes v(x,y) \ar[d]^-{\psi_{x,x,y}} \\
& v(x,y)}
$$
The morphisms $(v,\psi)\to (v',\psi')$ are sets 
$\{\xymatrix@C=35pt{
v(x,y) \ar[r]|-{\,\varphi_{x,y}\,} &
v'(x,y)}\}_{x,y\in X}$ 
of morphisms in $V$ for which the following diagram commutes for all $x,y,z\in
X$. 
$$
\xymatrix{
a(x,y)\otimes v(y,z) \ar[r]^-{1\otimes \varphi_{y,z}} \ar[d]_-{\psi_{x,y,z}}&
a(x,y)\otimes v'(y,z)\ar[d]^-{\psi'_{x,y,z}} \\
v(x,z) \ar[r]_-{\varphi_{x,z}} &
v'(x,z)}
$$
By \cite[Proposition 4.2]{BCV} this is a monoidal category with the product
$(v\otimes v')(x,y):=v(x,y)\otimes v'(x,y)$ for all $x,y\in X$ and
$$
\xymatrix@C=25pt@R=10pt{
a(x,y) \!\otimes\! (v\!\otimes\! v')(y,z) 
\ar[r]^-{\raisebox{10pt}{${}_{\delta_{x,y} \otimes 1}$}} &
a(x,y) \!\otimes\! a(x,y) \!\otimes\! v(y,z)\!\otimes\! v'(y,z) 
\ar[r]^-{\raisebox{10pt}{${}_{1\otimes c \otimes 1}$}} &
{\color{white} (v\!\otimes\! v')(x,z)}\\
& a(x,y)  \!\otimes\! v(y,z) \!\otimes\! a(x,y) \!\otimes\! v'(y,z) 
\ar[r]^-{\raisebox{10pt}{${}_{\psi_{x,y,z} \otimes \psi'_{x,y,z}}$}} &
(v\!\otimes\! v')(x,z)}
$$
for $x,y,z\in X$.
\item[{---}]
The monoidal full subcategory of the Eilenberg--Moore category of the
opmonoidal monad $\mathsf{Span}\vert V(X,a)$ on $\mathsf{Span}\vert V(X,X)$,
whose objects live on the complete $X$-span. 
\end{itemize}

\subsection{Hopf monads in $\mathsf{Span}\vert V$ versus Hopf categories}
Consider again a  $\mathsf{Cmd}(V)$-enriched category with
object set $X$ and hom objects $(a(x,y),\delta_{x,y},\varepsilon_{x,y})$ for
$(x,y)\in X \times X$. Denote the composition compatibility morphisms by
$\mu$ and denote the unit compatibility morphisms by $\eta$ as in the previous
section.
As we saw in Section \ref{sec:Span|V-bimonad}, it can be regarded equivalently
as a bimonoid in the duoidal category $\mathsf{Span}\vert V(X,X)$. In the
current situation the {\em antipode} in the sense of \cite[Theorem 7.2]{BoLa}
turns out to be a set of morphisms in $V$, 
$\{\xymatrix@C=30pt{
a(v,w) \ar[r]|(.47){\,\sigma_{v,w}\,} & a(w,v)}
\}_{v,w\in X}$, 
subject to the axioms in \cite[Theorem 7.2]{BoLa}. The first antipode axiom in
\cite[Theorem 7.2]{BoLa} takes now the form in Figure \ref{fig:antipode}. In
that figure, for natural numbers $n\geq m$, we denote by $p_m$ the
$m^{\textrm{th}}$ projection from the $n$-fold Cartesian product of $X$ to
$X$, sending $(q_1,\dots,q_n)$ to $q_m$.      
\begin{figure} 
\centering
\begin{sideways}
\scalebox{.92}{
$\xymatrix@C=23pt@R=45pt{
\mbox{$\begin{array}{l}
(X \stackrel{p_1}\longleftarrow X\times X \stackrel{p_1}\longrightarrow X \\
(v,w) \mapsto a(v,w)) 
\end{array}$} 
\ar[rr]^-{(1,(v,w) \mapsto \varepsilon_{v,w})} 
\ar[d]_-{(1,(v,w)\mapsto \delta_{v,w})} &&
\mbox{$\begin{array}{l}
(X \stackrel{p_1}\longleftarrow X\times X \stackrel{p_1}\longrightarrow X \\
(v,w) \mapsto K ) 
\end{array}$} 
\ar[r]^-{(p_1,1)} &
\mbox{$\begin{array}{l}
(X = X= X \\
(v,w) \mapsto K ) 
\end{array}$} 
\ar[d]^-{(\Delta,v\mapsto \eta_v)} 
\\
\mbox{$\begin{array}{l}
(X \stackrel{p_1}\longleftarrow X\times X \stackrel{p_1}\longrightarrow X \\
(v,w) \mapsto a(v,w) \otimes a(v,w) )
\end{array}$}  
\ar[r]_-{\raisebox{-20pt}{${}_{(\Delta \times 1,1)}$}} & 
\mbox{$\begin{array}{l}
(X \stackrel{p_1}\longleftarrow X\times X\times X \stackrel{p_2}
\longrightarrow X \\
(v,z,w) \mapsto a(v,w) \otimes a(z,w) )
\end{array}$} 
\ar[r]_-{\raisebox{-20pt}{${}_{((v,z,w)\mapsto (v,w,z),(v,z,w)\mapsto 1
      \otimes \sigma_{z,w})}$}} &  
\mbox{$\begin{array}{l}
(X \stackrel{p_1}\longleftarrow X\times X\times X \stackrel{p_3}
\longrightarrow X \\
(v,w,z) \mapsto a(v,w) \otimes a(w,z) ) 
\end{array}$} 
\ar[r]_-{\raisebox{-20pt}{${}_{((v,w,z)\mapsto (v,z),(v,w,z) \mapsto \mu_{v,w,z})}$}} &
\mbox{$\begin{array}{l}
(X \stackrel{p_1}\longleftarrow X\times X \stackrel{p_2}\longrightarrow X \\
(v,w) \mapsto a(v,w) )
\end{array}$} 
}$
}
\end{sideways}
\caption{The first antipode axiom}
\label{fig:antipode}
\end{figure}

The second antipode axiom is handled symmetrically. Comparing these diagrams
with \cite[Definition 3.3]{BCV} we conclude by \cite[Theorem 7.2]{BoLa} that
for any braided monoidal category $V$, the following notions coincide. 
\begin{itemize}
\item[{---}] A {\em Hopf $V$-category} in \cite[Definition
  3.3]{BCV}. Explicitly, this means a $\mathsf{Cmd}(V)$-enriched category with
  some object set $X$ and hom objects
  $(a(p,q),\delta_{p,q},\varepsilon_{p,q})$ for $(p,q)\in X \times X$, 
  composition compatibility morphisms $\mu_{p,q,r}:a(p,r) \otimes a(q,r)\to
  a(p,q)$ and unit compatibility morphisms $\eta_p:K \to a(p,p)$,
  for all $p,q,r\in X$; equipped with a further set 
$\{\xymatrix@C=25pt{
a(p,q) \ar[r]|(.47){\,\sigma_{p,q}\,} &
a(q,p)}\}_{p,q\in X}$
of morphisms in $V$ rendering commutative the following diagrams for
all $p,q\in X$. 
$$
\xymatrix{
a(p,q) \ar[r]^-{\delta_{p,q}} \ar[dd]_-{\varepsilon_{p,q}} &
a(p,q) \otimes a(p,q) \ar[d]^-{1\otimes \sigma_{p,q}} \\
& a(p,q) \otimes a(q,p) \ar[d]^-{\mu_{p,q,p}} \\
K \ar[r]_-{\eta_p} &
a(p,p)}\qquad 
\xymatrix{
a(p,q) \ar[r]^-{\delta_{p,q}} \ar[dd]_-{\varepsilon_{p,q}} &
a(p,q) \otimes a(p,q) \ar[d]^-{\sigma_{p,q}\otimes 1} \\
& a(q,p) \otimes a(p,q) \ar[d]^-{\mu_{q,p,q}} \\
K \ar[r]_-{\eta_q} &
a(q,q)}
$$
\item[{---}] A Hopf monad in $\mathsf{Span}|V$ on the naturally Frobenius 
opmap monoidale $X$ of Section \ref{sec:Span|V-monoidale}, in whose 1-cell
part $X\to X$ the underlying span is the complete span 
$\xymatrix@C=10pt{
X & \ar[l] X\times X \ar[r] &X}$.
\end{itemize}

\subsection{The functorial relation of Hopf group monoids and Hopf 
categories to Hopf polyads}
Regarding a braided monoidal category as a monoidal bicategory 
with a single 0-cell, there is a monoidal pseudofunctor $V\to
\mathsf{Cat}$ as follows. 

The single 0-cell of the bicategory $V$ is sent to the category $V$. 
A 2-cell in the bicategory $V$ --- that is, a morphism $f:p\to q$
in the category $V$ --- is sent to the natural transformation 
$f \otimes (-):p \otimes (-) \to q \otimes (-)$ between endofunctors on $V$.
This is clearly a pseudofunctor. It is monoidal as well via the following
ingredients. The unit-compatibility pseudo natural transformation is provided
by the 1-cell of $\mathsf{Cat}$ (i.e. functor) from the terminal category to
$V$ sending the only object to the monoidal unit $K$; and the isomorphism
$K\otimes K \cong K$ in $V$. The product-compatibility pseudonatural 
transformation has the object part provided by the monoidal product 
$\otimes: V \times V \to V$ and the morphism part given by the 
braiding $c$ of $V$ as $1\otimes c \otimes 1: p \otimes (-) \otimes q 
\otimes (-) \to p \otimes q \otimes (-) \otimes (-)$ for any object $(p,q)$ of
$V\times V$. The associativity and unitality modifications are induced by the
associativity and unitality natural isomorphisms of $V$. 

This monoidal pseudofunctor $V\to \mathsf{Cat}$ induces a monoidal
pseudofunctor from $\mathsf{Span}|V$ to $\mathsf{Span}|\mathsf{Cat}$ whose
unit- and product-compatibilities are pseudonatural transformations as
well. Since such monoidal pseudofunctors preserve monoidales (but not
necessarily opmap mono\-idales!), monads and opmonoidal morphisms, as well as
the invertibility of 2-cells, we conclude that they preserve Hopf monads. In
particular, the above monoidal pseudofunctor $\mathsf{Span}|V \to
\mathsf{Span}|\mathsf{Cat}$ takes both Hopf group monoids and Hopf categories
to Hopf polyads. 
Hopf polyads in the range of this monoidal pseudofunctor
$\mathsf{Span}|V \to \mathsf{Span}|\mathsf{Cat}$ were termed {\em
  representable} in \cite[Section 7.2]{Bru}.


\bibliographystyle{plain}

\begin{thebibliography}{10}

\bibitem{AgMa}  
Marcelo Aguiar and Swapneel Mahajan, 
{\em Monoidal Functors, Species and Hopf Algebras.} 
CRM Monograph Series 29, American Math. Soc. Providence, 2010. 
Electronically available at: 
\href{http://www.math.tamu.edu/~maguiar/a.pdf}{http://www.math.tamu.edu/~maguiar/a.pdf}.

\bibitem{BCV}
Eliezer Batista, Stefaan Caenepeel and Joost Vercruysse, 
{\em Hopf Categories,} 
Algebr. Represent. Theory 19 no. 5 (2016), 1173-1216.

\bibitem{Borceaux}
Francis Borceux, 
{\em Handbook of Categorical Algebra: Volume 1, Basic Category Theory.}
Cambridge University Press 1994.

\bibitem{BoLa}
Gabriella B\"ohm and Stephen Lack,
{\em Hopf comonads on naturally Frobenius map-monoidales,} 
J. Pure Appl. Algebra 220 no. 6 (2016), 2177-2213. 

\bibitem{BNSz} 
Gabriella B\"ohm, Florian Nill and Korn\'el Szlach\'anyi,
{\em Weak Hopf algebras. I. Integral theory and $C^*$-structure,} 
J. Algebra, 221 no. 2 (1999) 385-438. 

\bibitem{Bru}
Alain Brugui\`eres,
{\em Hopf polyads,}
preprint available at 
\href{http://arxiv.org/abs/1511.04639}{http://arxiv.org/abs/1511.04639}.

\bibitem{BLV} 
Alain Brugui\`eres, Steve Lack and Alexis Virelizier, 
{\em Hopf monads on monoidal categories,} 
Adv. Math. 227 no. 2 (2011) 745-800. 

\bibitem{BV}
Alain Brugui\`eres and Alexis Virelizier, 
{\em Hopf monads,} 
Adv. Math. 215 no. 2 (2007) 679-733.

\bibitem{CL}
Stefaan Caenepeel and Michel De Lombaerde, 
{\em A categorical approach to Turaev's Hopf group-coalgebras,} 
Comm. Algebra 34 (2006), 2631-2657.  

\bibitem{CLS}
Dimitri Chikhladze, Stephen Lack and Ross Street, 
{\em Hopf monoidal comonads,} 
Theory Appl. Categ. 24 no. 19 (2010) 554-563. 

\bibitem{GPS}
Robert Gordon, Anthony John Power and Ross Street, 
{\em Coherence for tricategories.} 
Mem. Amer. Math. Soc. 117 no. 558, 1995.

\bibitem{Gray}
John W. Gray,
{\em Formal Category Theory: Adjointness for 2-Categories.}
Springer LNM 391, 1974.

\bibitem{Lop-Fr}
Ignacio L\'opez Franco, Ross Street and Richard J. Wood, 
{\em Duals invert,} 
Appl. Categ. Structures 19 no. 1 (2011) 321-361.

\bibitem{Lop-FrStreetWood}
Ignacio L\'opez Franco, 
{\em Formal Hopf algebra theory. I. Hopf modules for pseudomonoids,} 
J. Pure Appl. Algebra 213 no. 6 (2009) 1046-1063.

\bibitem{McCr}
Paddy McCrudden,
{\em Opmonoidal monads.}
Theory Appl. Categ. 10 no 19 (2002) 469-485.

\bibitem{Moerdijk}
Ieke Moerdijk,
{\em Monads on tensor categories.}
J. Pure Appl. Algebra 168 no. 2–3 (2002) 189-208.

\bibitem{Rav}  
Douglas C. Ravenel, 
{\em Complex cobordism and stable homotopy groups of spheres,} 
Pure Appl. Math., vol. 121, Academic Press, Orlando, FL, 1986.

\bibitem{Sch}
Peter Schauenburg, 
{\em Duals and doubles of quantum groupoids ($\times_R$-Hopf algebras),} 
AMS Contemp. Math. 267 pp 273-299, AMS, Providence, RI, 2000.

\bibitem{Str:FTMI}
Ross Street, 
{\em The formal theory of monads,}
J. Pure Appl. Algebra 2 no. 2 (1972) 149-168.

\bibitem{St}
Ross Street, 
{\em Monoidal categories in, and linking, geometry and algebra,} 
Bull. Belg. Math. Soc. Simon Stevin, 19 no. 5 (2012) 769-821. 

\bibitem{T}
Vladimir G. Turaev, 
{\em Homotopy field theory in dimension 3 and crossed group-categories,} 
preprint available at 
\href{http://arxiv.org/abs/0005291}{http://arxiv.org/abs/0005291}.

\bibitem{Z}
Marco Zunino, 
{\em Double construction for crossed Hopf coalgebras,} 
J. Algebra 278 (2004), 43-75.

\end{thebibliography}

\end{document}